\documentclass[11pt,reqno]{amsart}
\usepackage[foot]{amsaddr}
\usepackage{import}
\usepackage[a4paper,margin=1.0in]{geometry}
\usepackage{bm}
\usepackage{psfrag}
\usepackage[usenames,dvipsnames]{color}
\usepackage{xcolor}
\usepackage[all,cmtip,line]{xy}
\usepackage{tikz} 
\usepackage[utf8]{inputenc}
\usepackage{pgfplots}
\usepackage{subfigure}
\usepackage{amssymb}
\usepackage{caption}
\usepackage{amsmath}
\usepackage[shortlabels]{enumitem}
\usepackage{comment}
\usepackage{algorithm}
\usepackage{algpseudocode}
\usepackage{mathrsfs}

\usepackage{mathtools}
\mathtoolsset{showonlyrefs}
\numberwithin{equation}{section}

\usepackage[T3,T1]{fontenc}
\DeclareSymbolFont{tipa}{T3}{cmr}{m}{n}
\DeclareMathAccent{\invbreve}{\mathalpha}{tipa}{16}
\DeclareMathOperator*{\argmin}{arg\,min}

\newtheorem{theorem}{Theorem}[section]
\newtheorem{lemma}[theorem]{Lemma}
\newtheorem{corollary}[theorem]{Corollary}
\newtheorem{proposition}[theorem]{Proposition} 

\newtheorem{problem}[theorem]{Problem}
\newtheorem{definition}[theorem]{Definition}

\newtheorem{remark}[theorem]{Remark}
\newtheorem{assumption}[theorem]{Assumption}

\newcommand{\sU}{{\normalfont\text{U}}}
\newcommand{\x}{{\bf{x}}}
\newcommand{\y}{{\bf{y}}}

\renewcommand{\div}{\operatorname{div}}

\definecolor{forestgreen}{rgb}{0.13, 0.55, 0.13}

\newcommand{\jp}[1]{{\color{blue}#1}}

\newcommand{\todo}[1]{{\color{red}[#1]}}

\newcommand{\fh}[1]{{\color{forestgreen}{[#1]}}}

\newcommand{\be}{\begin{equation}}
\newcommand{\ee}{\end{equation}}

\usepackage{algorithm}
\usepackage{algpseudocode}

\newcommand{\ba} {\bm a}

\newcommand{\bmu} {\bm{\mu}}

\newcommand{\N}{\mathrm{N}}

\newcommand{\bg}{\bm{g}}

\newcommand{\dual}[2]{\left\langle#1,#2\right\rangle}

\newcommand{\norm}[1]{\left \lVert #1 \right \rVert}
\newcommand{\snorm}[1]{\left \lvert #1 \right \rvert}

\renewcommand{\div}{{\rm div}}

\newcommand{\IC}{{\mathbb C}}

\newcommand{\IN}{{\mathbb N}}

\newcommand{\IR}{{\mathbb R}}
\newcommand{\IU}{{\mathbb U}}
\newcommand{\IT}{{\mathbb T}}
\newcommand{\IW}{{\mathbb W}}

\newcommand{\IY}{{\mathbb Y}}

\newcommand{\bp}{{\bm p}}

\newcommand{\fpcs}[1]{\operatorname{\mathsf{f}\pc}}

\newcommand{\cmspaceh}[4]{\mathcal{C}^{#1,#2} \left( #3, #4 \right)}

\newcommand{\rgeo}[1]{\mathcal{C}_b^{#1}\left( (-1,1), \IR^2 \right)}

\newcommand{\rgeoh}[2]{\mathcal{C}_b^{#1,#2}\left( (-1,1), \IR^2 \right)}

\renewcommand{\div}{\operatorname{div}}

\newcommand{\cP}{\mathcal{P}}

\newcommand{\pc}{}

\newcommand{\bla}{\boldsymbol \lambda}
\newcommand{\bphi}{\boldsymbol \phi}

\newcommand{\VH}{\bm{H}}

\newcommand{\bn}{\bm{n}}
\newcommand{\bu}{\bm{u}}

\newcommand{\bc}{\bm{c}}

\newcommand{\bz}{\bm{z}}

\newcommand{\bv}{\bm{v}}
\newcommand{\br}{\bm{r}}

\newcommand{\bx}{\bm{x}}
\newcommand{\by}{\bm{y}}

\newcommand{\href}{\VH^{\mathrm{ref}}}

\title[Reduced Basis for Multiple Open Arcs]{Reduced Basis Method for the Elastic Scattering by Multiple Shape-Parametric Open Arcs in Two Dimensions}

\author{Jos\'e Pinto$^\dagger$}
\address{$^\dagger$Facultad de Ingenier\'ia y Ciencias, Universidad Adolfo Ib\'a\~nez, Santiago, Chile}
\author{Fernando Henr\'iquez$^\ddagger$}
\address{$^\ddagger$Chair of Computational Mathematics and Simulation Science (MCSS), \'Ecole Polytechnique F\'ed\'eral de Lausanne, Lausanne, Switzerland.}
\thanks{This work was funded by the ANID grant Fondecyt Iniciaci\'on N°11230248.}
\email{jose.pinto@uai.cl}
\email{fernando.henriquez@epfl.ch}

\begin{document}
\maketitle

\begin{abstract}
We consider the elastic scattering problem by multiple disjoint arcs or \emph{cracks} in two spatial dimensions. A key aspect of our approach lies in the parametric description of each arc's shape, which is controlled by a potentially high-dimensional, possibly countably infinite, set of parameters. 
We are interested in the efficient approximation of the parameter-to-solution map employing model order reduction techniques, specifically the reduced basis method.

Initially, we utilize boundary potentials to transform the boundary value problem, originally posed in an unbounded domain, into a system of boundary integral equations set on the parametrically defined open arcs. Our aim is to construct a rapid surrogate for solving this problem. To achieve this, we adopt the two-phase paradigm of the reduced basis method. In the offline phase, we compute solutions for this problem under the assumption of complete decoupling among arcs for various shapes. Leveraging these high-fidelity solutions and Proper Orthogonal Decomposition (POD), we construct a reduced-order basis tailored to the single arc problem. Subsequently, in the online phase, when computing solutions for the multiple arc problem with a new parametric input, we utilize the aforementioned basis for each individual arc.
To expedite the offline phase, we employ a modified version of the Empirical Interpolation Method (EIM) to compute a precise and cost-effective affine representation of the interaction terms between arcs. Finally, we present a series of numerical experiments demonstrating the advantages of our proposed method in terms of both accuracy and computational efficiency.
\end{abstract}

%
%
%


\section{Introduction}
Solving parametric partial differential equations (pPDEs) is a common task in many areas of science and engineering.
Parameters are used to describe various aspects of a mathematical models, for example, material properties
or variations in the problem's physical domain of definition. Traditionally, methods such as Finite Differences, Finite Volumes,
or Finite Elements have been employed to numerically solve these problems.
However, applications such as multiple-query or real-time problems require a repeated
and swift evaluation of the high-fidelity or full-order model for different parametric inputs. 
This rapidly becomes computationallly untractable, therefore complexity reduction methods 
in the parameter space are required for a fast and efficient treatment of these computational models.

The Reduced Basis (RB) method aims at accelerating the computation
of the solution of pPDEs by using a two stages paradigm.
First, in an \emph{offline} phase a collection
of so-called \emph{snapshots} or {high-fidelity} solutions
of the problem are computed for a number of parametric inputs.
Then, a basis of reduced dimension is computed using this collection of solutions.
Currently, there are two main approaches of doing so: POD \cite{liang2002proper}
and greedy strategies \cite{hesthaven2014efficient,bui-thanh_model_2008,buffa_priori_2012,devore_greedy_2013}.
POD establishes in advance a set of samples in the parameter space to calculate the corresponding high-fidelity 
solutions. A reduced basis is then computed using the Singular Value Decomposition (SVD) of the snapshot matrix.
Conversely, greedy strategies involve adding a new element to be basis one after another in a sequential manner
in a such a way that new element provides the best improvement
in \emph{solution's manifold} approximation.
Once the reduced basis is constructed using either method, the original high-fidelity problem is projected onto this basis.
This step represents the \emph{online} phase of the RB method.
We direct the reader to \cite{hesthaven2016certified,quarteroni2015reduced,prud2002reliable,rozza2014fundamentals} for 
further details on these approaches. 

The construction of reduced spaces becomes a computational challenge particularly 
when the parameter space is of high, possibly countably infinite, dimension. Indeed,
the efficient approximation of maps with high-dimensional parametric inputs poses 
major challenges to traditional computational methods due to the so-called \emph{curse of dimensionality}
in the paramter space. In \cite{CCS15}, polynomial surrogates of high-dimensional
input maps are shown to converge independently of the dimension provided that there exists an \emph{holomorphic}
or \emph{analytic} dependence upon the parametric input.
Computationally, this property, usually referred to as parametric holomorphy, 
is the foundation to state provably dimension-independent convergence rates for a variety of methods including, for instance,
Smolyak interpolation and quadrature \cite{ZS17,SS13} and  
high-order Quasi-Monte Carlo integration (HoQMC) \cite{DKL14,DGGS16} among others. 
In particular, in the context of the model order reducton and the RB method, parametric holomorphy
gurantees dimension-inpedent bounds for the Kolmogorov's width \cite{cohen2016kolmogorov,CD15}, which in turn
imply dimension-independent convergence rates of both the POD-based and greedy strategies
for the RB method \cite{binev2011convergence,chen2015sparse,chen2016sparse,chen2016adaptive,buffa2012priori,maday2002priori}.

In the context of shape-parametric Boundary Integral Operators (BIOs), one may find a variety 
of works addressing and proving the aforementioned parametric
holomorphy property.
In \cite{Henriquez2021,henriquez2021_thesis} the holomorphic dependence
of the Calder\'on projector upon boundaries of class $\mathcal{C}^2$ has been established.
Furthermore, in \cite{dalla2022multi,DLM22,dalla2013family} analytic shape differentiability of Os
has been studied for problems in two and three dimensions.
In \cite{DH23_23_1} parametric holomorphy of the combined integral operator has been proved
for piece-wise $\mathcal{C}^2$-domains, thus allowing polygonal/polyhedral boundaries.
However, relevant to the current work is \cite{pinto2024shape}, in which
the case of multiple open arcs in two-dimensions has been addressed. 

\subsection{Contributions}
In this work, we perform model order reduction for the elastic scattering problem by multiple shape-parametric arcs.
Firstly, as in \cite{stephane,stephan1984augmented,sloan1991,JHP20,kress1996,kress2000}, 
we cast the original boundary value problem as a system of boundary integral equations (BIEs) posed on the collection of open arcs.
Then, following the approach presented in \cite{GHS12}, we propose and thoroughly analyze
a reduced basis method for said shape-parametric formulation.
A key insight of this method consists in the construction 
of a reduced basis for each shape-parametric arc, which is then used as a building block for the complexity
reduction of the multiple interacting arcs configuration.
For the numerical approximation of the high-fidelity solution,
we use a custom spectral Galerkin Boundary Element (BE) implementation 
which is tailored to deal with the problem's characteristic singular nature at the arcs' endpoints.
We provide a comprehensive convergence analysis accounting for
the discretization in the parameter space, the BE discretization, 
and parametric dimension truncation.
Unlike previous works in this field, we also provide a systematic 
analysis for the construction of the RB space for the elastic scattering by multiple arcs
by using as a building block the single arc problem. These insights 
are supported by a series of numerical experiments.

\subsection{Outline}
This work is structured as follows. In Section \ref{sec:prelimnaries}, we introduce
notation and relevant technical results to be used throughout this work. 
Section \ref{sec:model_problem} introduces the elastic scattering problem
by multiple open arcs, together with its boundary integral formulation and
details of the spectral Galerkin BE discretization.
Next, in Section \ref{eq:rom_multiple_arcs} we introduce the reduced basis method for pPDES.
Particular emphasis is given to the construction of a reduced order model for the multiple arc
problem by taking as starting point a reduced basis constructed for a single arc.
In Section \ref{sec:convergence_analysis} we provide a complete analysis of the reduced
method for the multiple arc problem.
Finally, in Section \ref{sec:numerical_results} we present numerical experiments portraying
the performance of the reduced basis approach, whereas in Section \ref{sec:conclu}
we provide final concluding remarks.

\section{Preliminaries}\label{sec:prelimnaries}

\subsection{Notation}
Throughout, vectors are indicated by boldface symbols.
For any $\bv \in \IC^n$, with $n \in \IN$, we consider
the Euclidean norm $\norm{\bv} = \sqrt{\bv \cdot \overline{\bv}}$.
In $\IR^2$, we denote by $\bm{e}_1 = (1,0)^\top$, $\bm{e}_2 = (0,1)^\top$
the canonical vectors.
Given an angle $\theta \in [0,2\pi)$, the corresponding directional vector is
$\bm{e}_\theta = (\cos \theta,\sin \theta) ^\top$.
The rotation matrix associated with the angle $\theta$ is 
\begin{align}\label{eq:rotmatrix}
	\bm{R}_\theta 
	= 
	\begin{pmatrix}
    	\cos\theta & - \sin \theta  \\
    	\sin \theta & \cos \theta
	\end{pmatrix},
\end{align}
thus we have that $\bm{e}_\theta =  \bm{R}_\theta \bm{e}_1$.

Given real numbers $a,b$, we say that $a \lesssim b $
if there exists a positive constant $c$, independent of the
variables relevant to the corresponding analysis,
such that $a \leq c b$. If $a \lesssim b$ and
$b \lesssim a$ we write $a \cong b$.

Let $X$ be a Banach space. Its anti-dual is denoted by $X^\star$, 
and the evaluation of an element $f \in X^\star$ on an element $x \in X$ is denoted by $\langle f, x \rangle$.
Given another Banach space $Y$, we denote by $\mathcal{L}(X,Y)$ the space of bounded linear operators
from $X$ to $Y$. As it is customary, we equip it with the standard operator norm,
thus rendering it a Banach space itself.

Given a Hilbert space $H$ and a closed subspace $V \subset H$, we denote the corresponding orthogonal projection by $P^H_V$, or simply $P_V$ when the space $H$ is clear from the context.

\subsection{Functional Spaces} 
First, let us recall the definition of H\"older spaces.
Given tho non-empty, opened and connected sets
$\Omega_1 \subset \IR^{d_1}, \Omega_2 \subset \IR^{d_2}$, for $d_1, d_2 \in \IN$
the space $\cmspaceh{m}{\alpha}{\Omega_1}{\Omega_2}$ consists
of functions $f:  \Omega_1 \rightarrow \Omega_2$ with derivatives up
to order $m$ in $\Omega_1$, each of them continuous over $\overline{\Omega_1}$, 
and such that the derivatives of order $m$ fullfil the $\alpha$-H{\"o}lder continuous condition,
i.e.~$\|g(\x) -g(\y) \| \lesssim \|\x-\y\|^\alpha$ with $g$ being any derivative of order $m$ of $f$. 

In order to describe the relevant geometries we will make use of a subset of
$\cmspaceh{m}{\alpha}{(-1,1)}{\IR^2}$, denoted by $\rgeo{m,\alpha}$
consisting of all $\br \in \cmspaceh{m}{\alpha}{(-1,1)}{\IR^2}$
such that $\|\br'(t)\|>0$, $t\in (-1,1)$, and having a globally defined inverse. 

We also introduce the functional spaces used to properly state the elastic scattering problem on cracks.
Set $w(t) = \sqrt{1-t^2}$, $t\in(-1,1)$. We denote by $T_n(t)$ the $n$--th Chebyshev polynomial of
the first kind normalized according to 
\begin{equation}
	\int_{-1}^1 T_m(t) T_l(t) w^{-1}(t) dt 
	= 
	\delta_{m,l}, \quad l,m \in \IN_0,
\end{equation}
where $\delta_{l,m} = 1$ if $l=m$ and $\delta_{l,m}= 0$ if $l \neq m$.
For a smooth function $u: [-1,1] \rightarrow \IC$ we define two sequences
of Chebyshev coefficients as
\begin{equation}
	u_n 
	\coloneqq
	\int_{1}^1 u(t) T_n(t) dt
	\quad
	\text{and}
	\quad 
	\widehat{u}_n 
	\coloneqq 
	\int_{-1}^1 u(t) T_n(t) w^{-1}(t) dt, \quad n \in \IN_0.
\end{equation}
By using a duality argument, these definitions are extended to distributions,
and we define the following spaces: For $s \in \mathbb{R}$ we set
\begin{align}
	T^s 
	\coloneqq
	\left\lbrace 
		u : 
		\norm{u}_{T^s}^2 
		= 
		\sum_{n=0}^\infty (1+n^2)^s \snorm{{u}_n}^2 
		< 
		\infty
	\right\rbrace,  \\
	W^s
	\coloneqq
	\left\lbrace 
		u : 
		\norm{u}_{W^s}^2 
		= 
		\sum_{n=0}^\infty (1+n^2)^s \snorm{\widehat{u}_n}^2 
		< 
		\infty
	\right\rbrace.
\end{align}
For any $s \in \IR$, the dual space of $T^s$ can be identified with $W^{-s}$,
where the duality product is the extension of the $L^2((-1,1))$ inner product,
and throughout this work we assumed that this identification has been made. 

For certain values of $s$, these space coincide with standard Sobolev spaces
in the interval $(-1,1)$. In fact, we have that
\begin{equation}
	T^{-\frac{1}{2}} 
	= 
	\widetilde{H}^{-\frac{1}{2}}((-1,1))
	\quad 
	\text{and} 
	\quad 
	W^{\frac{1}{2}} 
	= 
	H^{\frac{1}{2}}((-1,1))
\end{equation}

Finally, we define the following product spaces 
\begin{align*}
\IT^s  = T^s \times T^s	\quad 
	\text{and} 
	\quad \quad \IW^s = W^s \times W^s.
\end{align*}
These will be extensively used in the next sections.

\subsection{Problem Geometry}\label{sec:geo}
We provide a precise description of the type of open arcs to be considered in the rest of the work. 
We uniquely identify each open arc with one of its corresponding parametrizations.
Each open arc is described as a function from $[-1,1]$ with values $\IR^2$
satisfying the following properties:
\begin{itemize}
	\item[(i)]
	The function is an element of $\mathcal{C}^{m,\alpha}((-1,1);\IR^2)$,
	with $m \in \IN$, $\alpha \in [0,1]$, and $m + \alpha > 2$.
    	\item[(ii)]
    	The derivative of the function, i.e.~the tangent vector,
	is nowhere null. Therefore, it belongs to $\rgeo{m,\alpha}$.
\end{itemize}
The latter requirement implies that the parametrization function is invertible. 
Of special interest will be collection of open arcs that are determined by a parametric input
$\by \in \text{U} =  \left[ -\frac{1}{2}, \frac{1}{2} \right]^\IN$ of the form: 
\begin{equation}\label{eq:paramarc}
	\br({\by},t) 
	= 
	\br_0(t)
	+ 
	\sum_{n=1}^\infty y_n \br_n(t) ,
\end{equation}
here $\br_0$ is an open arc, and in the following we refer to it as the reference arc.
In the following, for each $\by \in \text{U}$, when referring to an open arc as an element of $\mathcal{C}^{m,\alpha}((-1,1);\IR^2)$
we use the notaton $\br({\by})$. In contrast, when referring to a point in $\IR^2$ described by the arc's parametrization we use the
notation $\br({\by},t)$ for some $t \in [-1,1]$.
The set $\{\br_n\}_{n \in \IN}$ is a subset of $\mathcal{C}^{m,\alpha}((-1,1),\IR^2)$ and we will refer to it as the perturbation basis.
In order to ensure that for any value of the parameter $\by$ the parametrization $\br({\by}$
as in\eqref{eq:paramarc} is indeed an open arc in $\IR^2$, we work under the following 
assumptions.
\begin{assumption}\label{assump:b_decay}
\hfill
\begin{enumerate}
	\item
	The perturbation basis $\br_n$ is such that, $b_n = \|{ \br}_n\|_{\mathcal{C}^{m,\alpha}((-1,1);\IR^2)}$,
	$n\in \IN$, is a sequence in $\ell^p(\IN)$ for some $p\in (0,1)$.
    	\item 
    	There exists $\zeta\in(0,1)$ such that
	\begin{align}
		\sup_{t\in (-1,1)} 
		\sum_{n=1}^\infty 
		\|(\br_n)'(t)\| 
		\leq
		\zeta
		\inf_{t \in (-1,1) }
		\|(\br_0)'(t)\| .
	\end{align}
\end{enumerate}
\end{assumption}

Under these conditions for $\br_0$ and $\{b_n\}_{n \in \IN}$
we can ensure that $\br({\by})$ is an open arc for any $\by \in \text{U}$. 

As we are intereseted in the multiple arc problem,
we also consider $M$ parametrically defined collections
of open arcs, each of them of the form \eqref{eq:paramarc}.
We denote these families by $\br^1({\by^1}), \hdots,\br^M({\by^M})$,
each of them having their corresponding reference arc denoted by $\br^1_{0} ,\hdots \br^M_{0}$,
while the perturbations basis would be assumed to be the same for any of the $M$ collections. 

We further assume that fixed parts are line segments of the form 
$$
 	\br^j_{0}(t) 
	= 
	\bm{c}_j + \varrho_j (\cos \varphi_j, \sin \varphi_j)^\top t,
	\quad
	t\in (-1,1),
	\quad
	j=1,\dots,M.
$$
where $\bm{c}_j \in [-B,B] \times [-B,B]$, for some $B >0$, is the arc's center.
The length of the segment is $\varrho_j \in [r_\text{min}, r_\text{max}]$,
for given $0<r_\text{min}< r_\text{max}$, while $\varphi_j \in [0, \pi)$ defines the orientation of the segment.

We also assume that there exist $d_\text{min}, d_\text{max}$,
with $0<d_\text{min}< d_\text{max}$, such that $d_\text{min}\leq \|\bm{c}_k-\bm{c}_j\| \leq d_\text{max}$, for $k \neq j$,
and that 
$$
	d_\text{min} 
	> 
	2 \left(r_\text{max}+ \sum_{n=1}^\infty  \sup_{t \in [-1,1]}\| \br_n(t) \|\right),
$$
This last condition ensures that the arcs are pairwise disjoint.

\section{Elastic Wave Scattering by Multiple Open Arcs}
\label{sec:model_problem}
In this section, we introduce the elastic scattering problem by multiple open arcs, togethwe with its spectral BE Galerkin discretization 
using weighed polynomials.
We refer to \cite{Tao22} and references therein for more details on these aspects of the problem.   

\subsection{Problem Formulation}

Given $\kappa \in \IR$, $\theta \in [0,2\pi)$ we denote a scalar
plane wave with incidence angle $\theta$ and wavenumber $\kappa$ as 
$$
	g_{\theta,\kappa}(\bx) 
	= 
	\exp
	\left(
		\imath \kappa \bm{e}_\theta \cdot \bx 
	\right).
$$

Let us consider $M$ families of parameterized open arcs as in Section \ref{sec:geo},
whose image in $\IR^2$ are be denoted by $\Gamma(\by^1,\hdots,\by^M)$.
We fix a direction $\theta_0 \in [0,2\pi)$ that corresponds to that of the incoming plane-wave.
In addition, we also fix $\omega$ as the the problem's frequency, and we denote by $\lambda,\mu$ the Lam\'e parameters.
Then, for a given realization of the parameters $\by^1 ,\hdots, \by^M$, we seek $\bm{U}: \IR^2  \rightarrow \IC^2$ such that 
\begin{align}\label{eq:volprob}
	(\mu \Delta+ (\lambda +\mu)\nabla \nabla \cdot)\bm{U} + \omega^2 \bm{U}  
	&
	= 
	0
	\quad \text{in } \ \IR^2 \setminus \overline{\Gamma(\by^1,\hdots,\by^M)},\\
    	\bm{U} 
	&= 
	\bm{e}_{\theta_0} 
	g_{\theta_0 \kappa_{\mathsf{p}}}, 
	\quad \text{on }  \Gamma(\by^1,\hdots,\by^M)\\
	\label{eq:radcond}
	&
	+\text{Condition at infinity},
\end{align}
where $\kappa_{\mathsf{p}} =\sqrt{ \frac{\omega}{\lambda+2\mu}} $, and the radiation condition at infinity is the standard Kupradze's one. 

Using standard arguments of boundary integral formulations \cite[Chapter 7]{mclean2000strongly} we look for a solution of the form 
\begin{equation}
    \bm{U}(\x) = \sum_{j=1}^M \int_{-1}^1 	{\mathbf G}(\x,\br^j({\by^j},\tau)) \bm{u}^j(\by^1,\hdots, \by^M,\tau) d\tau  , \quad \x \in \IR^2 \setminus \overline{\Gamma(\by^1,\hdots,\by^M)}, 
    \label{eq:repformula}
\end{equation}
where $\mathbf{G}$\footnote{The fundamental solution depends on the parameters $\omega, \lambda, \mu$, however as these are fixed for each instance of the problem we do not incorporate them in the notation.} 
is the fundamental solution to the elastic wave operator,
and $\bm{u}^j(\by^1,\hdots, \by^M,\cdot)$ are unknown densities
defined in $(-1,1)$, for $j=1,\hdots,M$.
By imposing the boundary conditions on the representation formula stated in \eqref{eq:repformula}
we obtain the following system of boundary integral equations for the unknown
densities $\bm{u}^j(\by^1,\hdots, \by^M,\cdot)$, for $j=1,\hdots,M$:
\begin{align}\label{eq:biesrecast}
	\sum_{j=1}^M \bm{V}_{\br^k({\by^k}),\br^j({\by^j})}\bm{u}^j (\by^1,\hdots,\by^M,\cdot) (t) 
	= 
	\bm{g}_k(t), \quad k = 1, \hdots, M,
\end{align} 
where the weakly singular operator $\bm{V}_{\br,\bp}$ is defined,
for a pair of open arcs $\br, \bp$ and a density function $\bm{u}$, as 
\begin{equation*} \bm{V}_{\br,\bp} 
	\bm{u}(t)
	=
	\int\limits_{-1}^{1} \mathbf{G}(\br(t),\bp(\tau)) \bm{u}(\tau) d\tau 
\end{equation*}
and the right-hand side is given by 
\begin{equation}
\label{eq:rhsnot}
    \bm{g}_k(t) = \bm{e}_{\theta_0} g_{\theta_0, \kappa_p}(\br^k(\by^k)(t)).
\end{equation}

\subsection{High-Fidelity Discretization} 
\label{sec:high_fidelity_solver}
We now provide the construction of the spectral BE Galerkin discretization of \eqref{eq:biesrecast} used as
the high-fidelity solver in this work. 
To this end, we define the following family of finite dimensional spaces
\begin{align*}
	T_N 
	= \left\lbrace 
		u(t) = \sum_{n=0}^N a_n T_n(t) w^{-1}(t)\  \text{ with }\ \{a_n\}_{n=0}^N \subset \IC \right\rbrace, \quad N \in \IN,
\end{align*}
and set $\IT_N = T_N \times T_N$, for $N \in \IN$.
Having introducing these spaces, the discrete version of \eqref{eq:biesrecast} reads as follows:
For $\by^1,\hdots, \by^M \in \text{U}$, we seek
$u^j_N(\by^1,\hdots,\by^M, \cdot) \in \IT_N$, such that 
\begin{align}\label{eq:discbiesrecast}
	\sum_{j=1}^M P_{\IT_N}\bm{V}_{\br^k({\by^k}),\br^j({\by^j})}\bm{u}^j_N (\by^1,\hdots,\by^M,\cdot) (t) 
	= 
	P_{\IT_N}\bm{g}_k(t), \quad k = 1, \hdots, M.
\end{align} 
The following formulation presents the algebraic form of this last equation. 
\begin{problem}
\label{prob:linsistemN}
    Given $N \in \IN$ we seek $\bm{a}^{1,p}(\by^1,\hdots,\by^M),\hdots, \bm{a}^{M,p}(\by^1,\hdots,\by^M) \in \IC^{N+1}$, for $p \in \{1,2\}$, such that 
    \begin{align*}
        \sum_{j=1}^M\mathbb{A}^{p,q} _{k,j} \bm{a}^{j,q}(\by^1,\hdots,\by^M) 
	= 
	\bm{g}^p_{N,k}, 
	\quad k \in  \{1,\hdots,M\}, \quad q \in \{1,2\},
    \end{align*}
    where 
$$
	\left(
		\mathbb{A}^{p,q}_{k,j}
	\right)_{\ell,m} 
	= 
	\dual{
		\bm{V}_{\br^k({\by^k}),\br^j({\by^j})} \frac{T_m}{w} 
		\bm{e}_q
	}{
		\frac{T_\ell}{w} 
		\bm{e}_p
	}, 
	\quad 
	\ell,m \in \{0,\hdots,N\}, \ p,q \in \{1,2\}.
$$
and
$$
	\left(
		\boldsymbol{g}^p_{N,k}
	\right)_{\ell} 
	= 
	\dual{
		\bm{g}_{k}
	}{
		\frac{T_l}{w} \bm{e}_p 
	},
	\quad 
	\ell \in \{0,\hdots,N\}, \ p \in \{1,2\}.
$$
\end{problem}
The relation between the solution of Problem \ref{prob:linsistemN} and the solution of \eqref{eq:biesrecast} is the following:
$$
	\bm{u}^j_{N}(\by^1,\hdots,\by^M,t) 
	= 
	\sum_{m=0}^N \sum_{p=1}^2 
	\left(\boldsymbol{a}^{j,p}\left(\by^1,\hdots,\by^M\right)\right)_{m} T_m(t){w(t)}^{-1} \bm{e}_p.
$$
\begin{remark}\label{rmk:convergence_galerkin_bem}
In \cite{Tao22} it was proven that if for a given realization of the parameters $\by^1,\hdots, \by^M$ the resulting open arcs $\br^1(\by^1), \hdots, \br^M(\by^M)$ are analytic, then asymptotically in $N$ we obtain exponential convergence, i.e. there exist $\rho >1$ and $N_0 \in \IN$ such that
for any $N\geq N_0$ it holds
$$
	\sum_{j=1}^M
	\| \bm{u}^j(\by^1,\hdots, \by^M)-\bm{u}_N^j(\by^1,\hdots, \by^M) \|_{\IT^{-\frac{1}{2}}} 
	\lesssim \rho^{-N},
$$
where, in principle, the hidden constant depends on the parameters $\by^1,\hdots, \by^M \in \normalfont\text{U}$. 
\end{remark}

\subsubsection{Numerical Implementation}
\label{sec:numimphf}
We present an overview of some aspects concerning the computational implementation
of the prevously described spectral Galerkin BEM discretzation,
as these are important for the construction of the reduced model.
We again refer to \cite{Tao22,JHP20} for some more details and improvements to the basis implementation. 

For the implementation of the spectral Galerkin method we need to approximate two types of integrals: 
\begin{align}
     \dual{ \bm{V}_{\br^k({\by^k}),\br^j({\by^j})} \frac{T_m}{w} \bm{e}_q}{ \frac{T_l}{w} \bm{e}_p}
     \quad
     \text{and}
     \quad
     \dual{\bm{g}_k}{\frac{T_l}{w} \bm{e}_p}.
\end{align}
Let us start with the second type, which corresponds to right-hand-side.
Since $\bm{g}_k$ is a smooth function we can simply compute these integrals using FFT.
Concretely, we first construct a vector of evaluations, 
\begin{align}
\label{eq:rhsint}
    \bm{\mathfrak{g}_{k}} = (\bm{g}_k(x_j))_{j=0,\hdots,N_c-1}, \quad x_j = \cos \left( \frac{\pi (N_c-1-j)}{N_c-1} \right),
\end{align}
then we apply the discrete Fourier transform to this vector\footnote{The parameter $N_c$ is selected according to a certain tolerance and at worst grows linearly with $N$.}, and finally we perform an escalation as in \cite{JHP20}.
The computational cost is of order $O(N\log N)$ per arc. 
Notice that we can think of the integration method as a linear function acting on vector of evaluations of the right-hand side 
\begin{remark}
\label{rem:abstractint}
While we do not provide all the details concerning the computation of the integrals, (we again refer to \cite{JHP20} for a detailed discussion), in a more abstract setting the computation of the integrals of the form 
$$
	\dual{ f }{ \frac{T_l}{w} }, 
$$
where $f$ is known function, can be though as the application of a linear map $L$ to the vector $\bm{\mathfrak{f}} = f(x_j)_{j=0,\hdots,N_c-1}$, where the points $x_j$ are as in \eqref{eq:rhsint}.
\end{remark}

For the the computaton of the first type of terms in \eqref{eq:rhsint} we consider two separate cases.
Firstly, for smooth component of the cross interaction betweem arcs, i.e.~when $k \neq j$,
we use a tensorization of the method applied to the right-hand side, with a cost of $O(N^2 \log N)$ operations.

Finally, the self-interaction case, k=j, is treated as in \cite{Tao22},
can be reduced to the computation of a regular integral (this is done as in the cross interaction case), plus an integral of the form 
\begin{align}\label{eq:sinint}
	I_{l,m} 
	= 
	\int_{-1}^1 \int_{-1}^1 \log|t-\tau| J(t,\tau) \frac{T_m}{w}(\tau) \frac{T_l}{w}(t) d\tau dt, 
\end{align}
where $J$ is an analytic function. Using the exact same procedure used to treat the crossed interactions together with
the orthogonality properties of Chebyshev polynomials, we can construct the approximation
$$
	J(t,s) \approx \sum_{p=0}^Q \sum_{q=0}^Q j_{p,q} T_p(t) T_q(s),
$$
with a computational cost of $O(Q^2 \log Q)$.
On the other hand, from the well-known expansion
$$
	\log |t-s| = \sum_{n=0}^\infty d_n T_n(t) T_n(s),
$$
using traditional Chebyshev properties we have that
\begin{align}
\label{eq:sumdir}
 I_{l,m} \approx \sum_{n=0}^\infty d_n (j_{l+n,m+n}+j_{l+n,|m-n|}+j_{|l-n|,m+n}+j_{|l-n|,|m-n|}),   
\end{align}
where the error decays exponential with respect to $Q$, which is in turn proportional to $N$.
The evaluation of \eqref{eq:sumdir} has a computational cost of $O(N^3)$.
However, this could be reduced to $O(N^2 \log N ) $
by using the convolution properties of the discrete Fourier transform. 
The total cost of assembling the linear system of equations is of order $O(M^2 N^2 \log N)$. 

\begin{remark}
\label{rem:matints}
As in Remark \ref{rem:abstractint}, the algorithm for the computation of the matrix terms can be thought as linear function acting on evaluations of a kernel function. In the cross-interactions case the kernel function is the fundamental solution, while for the self-interactions we have two different linear functions, one acting on the evaluations of the function $J(t,s)$ (the same as in \eqref{eq:sinint}), and another one, corresponding to the regular part, which is the difference between the fundamental solution and $J(t,s)\log|t-s|$.
\end{remark}
\section{Model Order Reduction for Multiple Arcs}
\label{eq:rom_multiple_arcs}
The high-fidelity discretization introduced in Section \ref{sec:high_fidelity_solver} provides
a fast method to approximate the solution of the elastic scattering problem
for a fixed geometry configuration determined by the parameters $\by^1, \hdots, \by^M \in \text{U}$.
However, when aiming to solve this problem for a
large number of parametric inputs, 
a new strategy is required, in particular as the number of arcs increases. 
To this end, we adopt a model order reduction perspective and resort 
to the reduced basis method. 

This section is organized as follows: In Subsection \ref{sec:pod_rom},
we revisit the fundamentals of the Galerkin Proper Orthogonal Decomposition
(Galerkin-POD) approach, a well-established technique for constructing efficient reduced bases in the context of parametric problems. Subsequently, in Subsection \ref{sec:romconstruction}, we outline a method for constructing an efficient Reduced Order Model (ROM) for configurations with multiple arcs. This approach involves building individual ROMs for each single arc while disregarding interactions with other arcs.

To fully harness the benefits of the reduced basis approach 
one needs an efficient and fast way of constructing the high-fidelity problem
projected in the reduced space. Even though we use an affine-parametric
representation for each open arc, this does not translate into an affine 
decomposition of the underlying reduced problem.
Following ideas introduced in \cite{GHS12}  we apply the Empirical Interpolation Method (EIM)
as described in Section \ref{sec:eimconstruction} ahead.

\subsection{Galerkin-POD and Reduced Order Modelling}
\label{sec:pod_rom}
Let us briefly recall the Galerkin-POD method. Our presentation 
follows mainly \cite[Chapter 3]{hesthaven2016certified} and 
\cite[Chapter 6]{quarteroni2015reduced}.

For each $\by \in \text{U}$, we seek for $u(\by) \in X $ such that
\begin{align}\label{eq:variational_problem}
	\mathsf{a}
	\left(
		u(\by)
		,
		v
		;
		\by
	\right)
	=
	\mathsf{g}(v,\by),
	\quad
	\forall
	v\in X,
\end{align}
where $X$ is a Hilbert Space, for each $\by \in \mathbb{U}$,
$\mathsf{a}(\cdot,\cdot;\by): X\times X\rightarrow \IC$, and
$\mathsf{g}(\cdot;\by)\in X^\star$ denotes a parameter-dependent
sesquilinear form and an anti-linear functional acting on $X$.
We also define the solution manifold, i.e. the set all possible solutions
to \eqref{eq:variational_problem} as
\begin{equation}
	\mathcal{M}
	\coloneqq
	\left\{
		u(\by) \in X
		:
		\;
		\by \in \text{U}
	\right\}
	\subset X.
\end{equation}
In addition, let $\{X_N\}_{N \in \IN}$ be a family of finite-dimensional
subspace of $X$, each one of dimension $N$, and let 
$\{\varphi_1,\dots,\varphi_N\} \subset X$ be a suitable basis
of $X_N$, i.e. $X_N = \text{span}\{\varphi_1,\dots,\varphi_N\}$.
The Galerkin discretization of \eqref{eq:variational_problem} reads:
For each  $\by \in \text{U}$, we seek $u_N(\by) \in X_N$, such that 
\begin{align}\label{eq:variational_problem_discrete}
	\mathsf{a}
	\left(
		u_N(\by)
		,
		v
		;
		\by
	\right)
	=
	\mathsf{g}(v,\by),
	\quad
	\forall
	v\in X_N,
\end{align}
while we define the \emph{discrete} solution manifold as
\begin{align}
	\mathcal{M}_N
	\coloneqq
	\left\{
		u_N(\by) \in V_N
		:
		\;
		\by \in \text{U}
	\right\}.
\end{align} 
Equivalently, the Galerkin discretization of \eqref{eq:variational_problem}
can be formulated as follows: For each $\by \in \text{U}$ we seek 
${\bm a}_N(\by) \in \IC^N$ such that 
\begin{equation}	
	\mathbb{A}_N(\by) {\bm a}_N(\by) 
	= 
	{\bm g}_N(\by), 
\end{equation}
where for each $\by \in \text{U}$
\begin{equation}
	\left(
		\mathbb{ A}_N(\by)
	\right)_{i,j} 
	= 
	\mathsf{a} \left( \varphi_j, \varphi_i; \by\right)
	\quad
	i,j= 1,\hdots,N
	\quad
	\text{and}
	\quad
	({\bm g}_N(\by))_i = \mathsf{g}(\varphi_i;\by)
	\quad
	i=1,\hdots,N.
\end{equation}
One can readily observe that
$$
	u_N(\by) 
	= 
	\sum_{j=1}^N ({\bm a}_N(\by))_j \varphi_j.
$$
As we are working in a finite-dimensional space, it holds
$$
	\norm{
		\sum_{n=1}^N c_n\varphi_N 
	}_X^2  \cong
	 \sum_{n=1}^N 
	 \snorm{
	 	c_n
	}^2,
$$
with a constant depending on $X_N$
hence the norm $\norm{\cdot}_X$ in $X_N$ and the vector $2$--norm
of the coefficients $\{c_n\}_{n=1}^N \in \mathbb{C}^N$ are equivalent
with hidden constants depending of $N$.

For the construction of the reduced basis, we seek the subspace
$V^{\text{(rb)}}_R$ of $X_N$ of dimension $R<N$ solution to
the following minimization problem
\begin{equation}\label{eq:optimal_reduced_space}
	V^{\text{(rb)}}_R
	=
	\argmin_{\substack{Z_R \subset X_N \\ \text{dim}(Z_R)\leq R}} 
	\int\limits_{\text{U}} 
	\norm{
		u_N(\by) 
		- 
		\mathsf{P}_{Z_R} 
		u_N(\by)
	}^2_X 
	\text{d}\by,
\end{equation}
The above problem can be formulated using
the algebraic form of the Galerkin discretization.
The solution is given by a matrix $\mathbb{V}^{\text{(rb)}}_R \in \mathbb{C}^{N \times R}$ whose orthogonal 
columns span $V^{\text{(rb)}}_R$ and is obtained as the solution of the following problem: 
\begin{equation}\label{eq:minimization_problem}
	\mathbb{V}^{\text{(rb)}}_R
	=
	\argmin_{\mathbb{W} \in \mathscr{V}_R} 
	\int\limits_{\text{U}}
	\norm{
		{\bm a}(\by) 
		- 
		\mathbb{W}
		\mathbb{W}^\dagger
		{\bm a}(\by)
	}^2
	\,
	\text{d}{\by},
\end{equation}
where
$
	\mathscr{V}_R 
	= 
	\left\{
		\mathbb{W} \in \mathbb{C}^{N \times R}: \, 
		\mathbb{W}^\dagger \mathbb{W} 
		= 
		\mathbb{I}_R
	\right \}.
$

While the solution of this problem is known (see, e.g., \cite[Proposition 6.3]{quarteroni2015reduced}),
in practical implementations one considers a discretized version of the
high-dimensional integral in \eqref{eq:minimization_problem}.
Provided a dimension truncation $s \in \IN$, we consider an equal weights, $N_t$-points quadrature rule in $\text{U}^{s}$ with points
$\left \{\by_{1},\dots,\by_{N_t}\right \} \subset \text{U}^{s}$.
This yields the following approximation
of \eqref{eq:minimization_problem}:
\begin{equation}\label{eq:minimization_rom_discrete}
	\mathbb{V}^{\text{(rb)}}_R
	=
	\argmin_{\mathbb{W} \in \mathscr{V}_R} 
	\frac{1}{N_t}
	\sum_{j=1}^{N_t}
	\norm{
		{\bm a}
		\left(
			\by_{j}
		\right) 
		- 
		\mathbb{W}
		\mathbb{W}^\dagger
		{\bm a}
		\left(
			\by_{j}
		\right)
	}^2
\end{equation}
The solution of this problem is given by the Schmidt-Eckart-Young
theorem (see, e.g., \cite[Proposition 6.1]{quarteroni2015reduced}).
Define the \emph{snapshot} matrix 
\begin{equation}
	\mathbb{S} 
	\coloneqq
	\begin{pmatrix}
		{\bm a}(\by_{1}) & \hdots & {\bm a}(\by_{N_t}) 
	\end{pmatrix}
	\in
	\mathbb{C}^{N \times N_t}
\end{equation}
and compute its SVD decomposition $\mathbb{S} = \mathbb{W} \Sigma \mathbb{Z}^\dagger$,
where
\begin{equation}
	\mathbb{W} 
	=
	\begin{pmatrix}
	 	{\bm w}_1, & \dots & ,{\bm w}_N
	\end{pmatrix}
	\in 
	\mathbb{C}^{N \times N}
	\quad
	\text{and}
	\quad
	\mathbb{Z} 
	=
	\begin{pmatrix}
	 	{\bm z}_1, & \dots & ,{\bm z}_{N_t}
	\end{pmatrix}
	\in 
	\mathbb{C}^{N_t \times N_t}
\end{equation}
are orthonormal matices.
Then $\mathbb{V}^{\text{(rb)}}_R$ in \eqref{eq:minimization_rom_discrete} 
is the matrix containing the $R$ vectors of $\mathbb{W}$ associated to the $R$ 
largest singular values of $\mathbb{S}$, which, under the assumption that the singular
values in $\Sigma$ are sorted in descreasing order.

Next, we define
\begin{equation}\label{eq:rb_space_V}
	\varphi^{\text{(rb)}}_i
	=
	\sum_{j=1}^{N}
	\left(
		{\bm w}_i
	\right)_j
	\varphi_j
	\quad
	\text{and}
	\quad
	V^{\text{(rb)}}_{R}
	=
	\text{span}
	\left\{
		\varphi^{\text{(rb)}}_1
		,
		\dots
		,
		\varphi^{\text{(rb)}}_R	
	\right\}.
\end{equation}
The Galerkin discretization 
of \eqref{eq:variational_problem} in the reduced space $V^{\text{(rb)}}_{R} \subset V_N$
reads: For each $\by \in \text{U}$, find $u_R^{(\text{rb})}(\by) \in V^{\text{(rb)}}_{R}$ such that 
\begin{align}\label{eq:variational_problem_discrete}
	\mathsf{a}
	\left(
		u^{\text{(rb)}}_R(\by)
		,
		v^{\text{(rb)}}_R
		;
		\by
	\right)
	=
	\mathsf{g}
	\left(
		v^{\text{(rb)}}_R
		;
		\by
	\right),
	\quad
	\forall
	v^{\text{(rb)}}_R \in V^{\text{(rb)}}_{R},
\end{align}
which in algebraic form reads
\begin{equation}
	\mathbb{A}^{\text{(rb)}}_R(\by)
	{\bm a}^{\text{(rb)}}_R(\by) 
	= 
	{\bm g}^{\text{(rb)}}_R(\by),
\end{equation}
where 
$$
	\left(
		\mathbb{A}^{\text{(rb)}}_R(\by)
	\right)_{\ell,m} 
	= 
	\mathsf{a}
	\left(
		\varphi^{\text{(rb)}}_m
		,
		\varphi^{\text{(rb)}}_\ell
		;
		\by
	\right)
	\quad
	\text{and}
	\quad
	\left(
		{\bm g}^{\text{(rb)}}_R(\by)
	\right)_\ell
	= 
	\mathsf{g}
	\left(
		\varphi^{\text{(rb)}}_\ell
		;
		\by
	\right), 
	\quad \ell,m = 1,\hdots,R,
$$
or equivalently 
\begin{align}\label{eq:naive_assmb}
	\mathbb{A}^{\text{(rb)}}_R(\by)
	= 
	\mathbb{V}^{\text{(rb)}\dagger}_R
	\mathbb{A}_N(\by)
	\mathbb{V}^{\text{(rb)}}_R
	\in 
	\mathbb{C}^{R \times R}
	\quad
	\text{and}
	\quad 
	{\bm g}^{\text{(rb)}}_R(\by) 
	= 
	\mathbb{V}^{\text{(rb)}\dagger}_R
	{\bm g}_N(\by).
\end{align}

 \begin{remark}[Criterium to select $R$]
 Provided a target tolerance $\epsilon_{\normalfont\text{svd}}$,
 we select $R$ as the smallest integer such that it holds
 $$
 	\frac{
 		\sum_{n=1}^R \sigma_n^2
 	}{
 		\sum_{n=R}^N \sigma_n^2
 	} 
 	>
 	1- \epsilon_{\normalfont\text{svd}}^2.
 $$
 where $\sigma_1\geq\cdots\geq \sigma_r>0$ are
 the singular values of $\mathbb{S}$, with $r  = \normalfont\text{rank}(\mathbb{S})$.
 \end{remark}

\subsection{Reduced Basis Construction}
\label{sec:romconstruction}

We now present how effectively apply the Galerkin POD
method discussed in Section \ref{sec:pod_rom} in the 
construction of a reduced basis for the multiple arc problem.  

Initially, we can attempt to directly apply the method to the
entire problem. However, its performance in addressing
the multiple arc problem is significantly hindered by two primary factors:
\begin{enumerate}
	\item 
	To construct the snapshot matrix, a significant number of parametric configurations
	with multiple arcs would be necessary.
	Each configuration is computationally demanding, especially as the number of arcs increases
    	\item 
   	The solution manifold (and its discrete counterpart)
	becomes difficult to approximate as it accounts 
	for shape variations of each crack.
\end{enumerate}
Even for a moderate number of arcs, a direct application
of the Galerkin-POD methods is prohibitively expensive.
Following \cite{GHS12}, we firstly construct a reduced basis
using the Galerkin POD method for a single arc.
Then, in the online stage, this basis is used as reduced order model for each individual arc.

\subsubsection{Reduced Basis Construction for a Single Arc}
\label{sec:reduced_basis_single_arc}
Using the notation of Section \ref{sec:geo}, we define a new family of parametrized open arcs as follows
\begin{equation}
\label{eq:gengeo}
	\bp(\by,t) 
	= 
	2B\begin{pmatrix}
    	y_1\\ y_2
	\end{pmatrix} 
	+ 	
	\varrho(y_3) 
	\begin{pmatrix}
	\cos \varphi(y_4) \\ \sin \varphi(y_4)
	\end{pmatrix} t 
	+ 
	\sum_{n\geq1}
	y_{n+4} 
	\br_n(t), \quad t \in [-1,1],
	\quad
	\by \in \text{U},
\end{equation}
where 
\begin{equation}
    \label{eq:geofunsaux}
    \varrho(z)= (r_\text{max}-r_\text{min})(z+\frac{1}{2})+r_\text{min}, \quad 
    \varphi(z)= \pi(z+\frac{1}{2}).
\end{equation}
We also define the corresponding set of all possible geometries 
and its dimension-truncated counterpart as 
\begin{equation}\label{eq:collection_of_arcs}
	\Sigma 
	=
	\left\{
		\bp(\by,\cdot): \; \by \in \text{U}
	\right\}
	\quad
	\text{and}
	\quad
	\Sigma_{s} 
	= 
	\left\{
		\bp(\by,\cdot): \;  \by \in \text{U}^{s}
	\right\},
\end{equation}
where as in the previous section we set $\text{U}^s = \left[ -\frac{1}{2}, \frac{1}{2} \right] ^s$

The construction of the reduced basis is as in Section \ref{sec:pod_rom},
however considering the following problems set on $\bp(\by)$:
For each $\by \in \text{U}$ we seek
$\bm{u}(\by),\widetilde{\bm{u}}(\by) \in \IT^{-\frac{1}{2}}$ such that
\begin{subequations}
\begin{align}
	\langle \mathbf{V}_{\bp(\by),\bp(\by)} \bm{u}(\by), \bm{v} \rangle 
	&
	= 
	\langle 
		\bm{e}_{\theta} g_{\theta,\kappa_\mathsf{p}} \circ \bp(\by)
		, 
		\bm{v}
	\rangle, 
	\quad \forall \bm{v} \in \IT^{-\frac{1}{2}},
	\label{eq:genbip1} \\
	\langle \mathbf{V}_{\bp(\by),\bp(\by)} \widetilde{\bm{u}}(\by), \bm{v} \rangle 
	&
	= 
	\langle \bm{e}_{\theta} g_{\theta+\frac{\pi}{2},\kappa_\mathsf{p}} \circ \bp(\by), \bm{v} \rangle, \quad \forall \bm{v} \in \IT^{-\frac{1}{2}},
	\label{eq:genbip2}
\end{align}
\end{subequations}
and its corresponding Galerkin discretizations using the method 
described in Section \ref{sec:high_fidelity_solver}. I.e.,
given $N \in \IN$, for each $\by \in \text{U}^{(s)}$
we seek $\bm{u}_N(\by),\widetilde{\bm{u}}_N(\by,\theta) \in \IT_N$ such that
\begin{subequations}
\begin{align}
	\dual{
		\mathbf{V}_{\bp(\by),\bp(\by)}
		\bm{u}_N(\by)
	}{
		\bm{v}_N
	}
	&
	= 
	\dual{
		\bm{e}_{\theta} g_{\theta,\kappa_\mathsf{p}} \circ \bp(\by)
	}{
		\bm{v}_N
	}, 
	\quad \forall \bm{v}_N \in \IT_{N}, 
	\label{eq:genbipd1}\\
	\dual{
		\mathbf{V}_{\bp(\by),\bp(\by)} \widetilde{\bm{u}}_N(\by)
	}{
		\bm{v}_N
	} 
	&
	= 
	\dual{
		\bm{e}_{\theta} g_{\theta+\frac{\pi}{2},\kappa_\mathsf{p}} \circ \bp(\by)
	}{
		\bm{v}_N
	}, 
	\quad \forall \bm{v}_N \in \IT_{N}.
	\label{eq:genbipd2}
\end{align}
\end{subequations}
For a given, fixed incident angle $\theta$, we consider the solution manifolds
\begin{equation}
	\mathcal{M}^\theta
	\coloneqq
	\left\{
	{\bm u}
	(\by)
	\in 
	\IT^{-\frac{1}{2}}:
	\;\;
	\by \in \text{U}
	\right\}
	\quad
	\mathcal{M}^{\theta+\frac{\pi}{2}}
	\coloneqq
	\left\{
	\widetilde{\bm u}
	(\by)
	\in 
	\IT^{-\frac{1}{2}}:
	\;\;
	\by \in \text{U}
	\right\}
\end{equation}
and its discrete counterparts
\begin{equation}
	\mathcal{M}^\theta_N
	\coloneqq
	\left\{
	{\bm u}_N
	(\by)
	\in 
	\IT_N:
	\;\;
	\by \in \text{U}
	\right\}
	\quad
	\mathcal{M}^{\theta+\frac{\pi}{2}}_N
	\coloneqq
	\left\{
	\widetilde{\bm u}_N
	(\by)
	\in 
	\IT_N:
	\;\;
	\by \in \text{U}
	\right\}.
\end{equation}
For the approximation of the discrete solution manifolds
$\mathcal{M}^\theta_N$ and $\mathcal{M}^{\theta+\frac{\pi}{2}}_N$,
we construct a snapshot matrix by sampling both problems in the parameter space. 
We then obtain a reduced space denoted by $V^{\text{(rb)}}_R$, with $R< 4(N+1)$. The reason as to why we include two different right-hand sides
is explained thoroughly ahead in Section \ref{sec:multiple_arcs_convergence}. 

\subsubsection{Reduced Basis for the Multiple Arc Problem}
\label{sec:reduced_basis_multiple_arc}
Let $V^{\text{(rb)}}_R$ be the reduced space of the single arc problem, constructed as in Section \ref{sec:reduced_basis_single_arc}.
For $\by^1,\hdots, \by^M \in \text{U}$, the reduced problem consists in seeking
$\bm{u}^{(\text{rb}),1}_R(\by^1,\hdots,\by^M)$, $\hdots$, $\bm{u}^{(\text{rb}),M}_R(\by^1,\hdots,\by^M) \in V^{\text{(rb)}}_R$ such that 
\begin{align}\label{eq:biesrecast}
	\sum_{j=1}^M 
	P^{\text{(rb)}}_R
		\bm{V}_{\br_k(\by_k),\br_j(\by_j)}
		{\bm u}^{(\text{rb}),j}_R (\by^1,\hdots, \by^M)	
	=  
	P^{\text{(rb)}}_R
	\bm{g}_k
	,
	\quad k = 1, \hdots, M. 
\end{align} 

Let $V^{\text{(rb)}}_{R} \subset V_N$ and
$\mathbb{V}^{\text{(rb)}}_R\in \mathbb{C}^{N\times R}$
be as in Section \ref{sec:reduced_basis_single_arc}.
The Galerkin problem in the reduced basis reads as follows.

\begin{problem}
\label{prob:rb}
    We seek $\bm{a}^{(\normalfont\text{rb}),1}(\by^1,\hdots,\by^M),\hdots, \bm{a}^{(\normalfont\text{rb}),M}(\by^1,\hdots,\by^M) \in \IC^{R}$, such that 
    \begin{align*}
        \sum_{j=1}^M\mathbb{A}^{(\normalfont\text{rb})}_R(\by^1, \hdots ,\by^M)_{k,j} \bm{a}^{(\normalfont\text{rb}),j}(\by^1,\hdots,\by^M) 
	= 
	\bm{g}^{(\normalfont\text{rb})}_{R,k}, 
	\quad k \in \{1,\hdots,M\},
    \end{align*}
    where 
$$
	\left(
		\mathbb{A}^{\normalfont(\text{rb})}_R(\by^1, \hdots ,\by^M)_{k,j}
	\right)_{\ell,m} 
	= 
	\dual{
		\bm{V}_{\br^k({\by^k}),\br^j({\by^j})} 
  \bm{\varphi}^{\normalfont\text{(rb)}}_m
	}{
		\bm{\varphi}^{\normalfont\text{(rb)}}_\ell
	}, 
	\quad 
	\ell,m \in \{1,\hdots,R\},
$$
and
$$
	\left(
		\bm{g}^{(\normalfont\text{rb})}_{R,k}
	\right)_{\ell} 
	= 
	\dual{
		\bm{g}_{k}
	}{
			\bm{\varphi}^{\text{\normalfont(rb)}}_\ell
	},
	\quad 
	\ell \in \{1,\hdots,R\}.
$$
The reduced basis $\bm{\varphi}^{\normalfont\text{(rb)}}_1,\hdots,\bm{\varphi}^{\normalfont\text{(rb)}}_R$
denotes the one constructed following the procedure of Section \ref{sec:pod_rom},
for the problem presented in Section \ref{sec:reduced_basis_single_arc}.
\end{problem}

%
%

\subsection{Reduced Basis Linear System Construction}
\label{sec:eimconstruction}
The newly constructed linear system, which stems from the
projection of the linear system of equations arising from the high-fidelity
model, is of substancial smaller size than that of the high-fidelity one.
However, to fully benefit from the reduced order basis method one
needs to be able to efficiently compute the full-order linear
system of equations in the online phase of the RB method.
In particular if the approach of \eqref{eq:naive_assmb} is used, the cost of the assembling the linear system for the high-fidelity model 
(see Section \ref{sec:high_fidelity_solver}) it would neglect any benefit of the reduced basis.

A common solution to improve the performance of the construction of the linear system, in the context where many evaluations for different parameters are required, is the Empirical Interpolation Method (EIM), see \cite[Chapter 5]{hesthaven2016certified}. In what follows we briefly introduce the latter method and explain how is used for the multiple arc problem.  

Given a function $f:(t,\by)\mapsto f(t,\by)$, where $t$ is called the physical variable, and $\by \in [-1,1]^{\IN}$ are the parameters representing the perturbations, the idea of the EIM is to construct an approximation of $f$ of the form: 
$$
f(t_j,\by) \approx \sum_{q=1}^Q c_q(\by) f(t_j,\by_q), \quad j=1,\hdots, N_c
$$
where $t_j,$ $j = 1,\hdots, N_c$  represent a set of pre-fixed points. The construction of this approximation is done in two stages. 

\begin{enumerate}
   	\item[(i)]
    	{\bf Offline Stage.}
    	Given a discretization of $\text{U}$, we find $Q \in \IN$,
	select the points $\by_q,$ $q=1,\hdots Q$ of the discretization
	of $\text{U}$, construct the functions $c_q(\cdot)$, $q=1,\hdots,Q$
	and store the values $f(t_j,\by_q)$ for future evaluations.
	This quantities are found using a greedy algorithm, in which
	is necessary to evaluate $f(t_j, \bz)$ for all the possible values
	of $t_j$, and $\bz$ in the discretization of $\text{U}$.
	While this can be very expensive, since the processes is
	independent of $\by$ (the parameter for which we want to
	evaluate the approximation), it has to be done only one time. 
    	\item[(ii)]
	{\bf Online Stage.}
	Given $\by \in \text{U}$ we evaluate the functions
	$c_q(\by)$, $q=1,\hdots,Q$ and the corresponding of
	approximations of $f(t_j,\by)$.
\end{enumerate}

Practical algorithms for the multiple arcs cases are given in the next section.
\subsubsection{EIM for Multiples Open Arcs}
We now explain how to use the EIM for the construction of the reduced
linear system for the multiple arcs problem. 

To solve the multiple arcs problem in the reduced basis space one needs to assemble the linear
system of equations described in Problem \ref{prob:rb}.
This implies the computation of $M^2$ block matrices: $M$ blocks accounting
for self-interaction terms, and $M^2-M$ for cross interaction between arcs.
It follows from \eqref{eq:naive_assmb} that each of these blocks are of the form
\begin{align}\label{eq:diraproach}
	\mathbb{A}_R^{(\text{rb})}(\by^1,\hdots,\by^M)_{k,j}
	=
	\mathbb{V}_R^{(\text{rb})\dagger}
		\mathbb{A}_N(\by^1,\hdots,\by^M)_{k,j}
	\mathbb{V}_R^{(\text{rb})},
\end{align}
where $\mathbb{A}_N(\by^1,\hdots,\by^M_{k,j}$ corresponds
to the matrix accounting for the interaction between arcs $k,j$ in the high-fidelity
space\footnote{This is the same matrix $\mathbb{A}_{k,j}$ of Section \ref{sec:high_fidelity_solver}, but we have made explicit the dependence of the parameters as its makes the use of the EIM more clear.}. Futhemore notice that while the full matrices $\mathbb{A}_R^{(\text{rb})}(\by^1,\hdots,\by^M)$, and $	\mathbb{A}_N(\by^1,\hdots,\by^M)_{k,j}$ depends on $M$ parameters, the block $k,j$ only depends on $\by^k$ and $\by^j$.
In the following, we simplify the notation but just including the two active parameters as arguments.

Let us consider first the computation of $\left(\mathbb{A}_R^{(\text{rb})}(\by^k,\by^j)\right)_{k,j}$, for
$k \neq j$. The the case $k=j$ is treated at the end of this section. 
According to Remark \ref{rem:matints} and \eqref{eq:diraproach}, there exists a linear 
function $L$, such that 
\begin{align}
	\mathbb{A}_N(\by^k,\by^j)_{k,j}
	\approx 
L(\bm{\mathfrak{G}}_{k,j}(\by^k,\by^j)),
\end{align}
where $\bm{\mathfrak{G}}_{k,j}(\by_1,\by_2)$, is a matrix constructed from
evaluations of the fundamental solution, i.e. consider $N_c$ points
$\{x_p\}_{p=1}^{N_c} \subset (-1,1)$
\begin{equation}
	\bm{\mathfrak{G}}_{k,j}(\by^k,\by^j) 
	= 
	\begin{pmatrix}
		\mathbf{G}(\br^k(\by^k)(x_1),\br^j(\by^j)(x_1)) & \cdots & \mathbf{G}(\br^k(\by^k)(x_1),\br^j(\by^j)(x_{N_c})) \\
		\vdots & \ddots & \vdots \\
		\mathbf{G}(\br^k(\by^k)(x_{N_c}),\br^j(\by^j)(x_1)) & \cdots & \mathbf{G}(\br^k(\by^k)(x_{N_c}),\br^j(\by^j)(x_{N_c})) 
	\end{pmatrix}.
\end{equation}
By introducing the linear function $L^R(\cdot) = (\mathbb{V}_R^{(\text{rb})})^\dagger L(\cdot)\mathbb{V}_R^{(\text{rb})}$, we can write an approximation of $\mathbb{A}_R^{(\text{rb})}(\by^k,\by^j)_{k,j}$, as
\begin{align}
	\mathbb{A}_R^{(\text{rb})}(\by^k,\by^j)_{k,j} 
	\approx 
L^R(\bm{\mathfrak{G}}_{k,j}(\by^k,\by^j)).
\end{align}
 Now we would use the EIM to form an  interpolation of $\bm{\mathfrak{G}_{k,j}}(\by_1,\by_2)$, this will result in an approximation of $  \mathbb{A}_R^{(\text{rb})}(\by^k,\by^j)_{k,j}$, once the linear map $L^R$ is applied to the interpolated evaluations. 

Without any further assumptions for each of the $M^2-M$ off-diagonal blocks
we have to compute a tailored interpolation based on the EIM for each pair of 
interactions. To reduce the 
computational burden, we instead form a global interpolation of  $\bm{\mathfrak{G}}_{k,j}(\by^k,\by^j)$ incorporating the local indices $k,j$ in the interpolation variables. 

To this end, given $s\in \IN$, and $\by^k,\by^j \in \text{U}^s$, we define $\bz^{k,j} \in \text{U}^{2s+6}$, whose components are 
\begin{align*}
	z^{k,j}_1 &=  \left( \frac{\varrho_k-r_\text{min}}{r_\text{max}-r_\text{min}}\right)-\frac{1}{2}\\
	z^{k,j}_2 &= \frac{ \varphi_k}{\pi} -\frac{1}{2}\\
	z^{k,j}_{n+2} &= (y^k)_{n}, \quad n = 1,\hdots,s\\
	z^{k,j}_{s+3} &=  \left( \frac{d-d_\text{min}}{d_\text{max}-d_\text{min}}\right)-\frac{1}{2}, \quad  d = \|\bm{c}_k-\bm{c}_j\|\\
	z^{k,j}_{s+4} &= \frac{ \theta}{\pi} -\frac{1}{2}, \quad \theta  = \arg(\bm{c}_k-\bm{c}_j)\\
	z^{k,j}_{s+5} &=  \left( \frac{\varrho_j-r_\text{min}}{r_\text{max}-r_\text{min}}\right)-\frac{1}{2}\\
	z^{k,j}_{s+6} &= \frac{ \varphi_j}{\pi} -\frac{1}{2}\\
	z^{k,j}_{n+s+6} &= (y^j)_n, \quad n= 1, \hdots,s
\end{align*}
where $r_\text{min},r_\text{max},d_\text{min},d_\text{max}$ are the global geometry parameters defined in Section \ref{sec:geo} and $\bm{c}_j,\bm{c}_k,\varrho_j,\varrho_k,\varphi_j,\varphi_k$ are the variables determining the fixed part of arcs $k,j$ respectively, also defined in the same section. For a general $\bz \in \text{U}^{2s+6}$ we introduce the following auxiliary parametrizations: 
\begin{subequations}
\begin{align}
	\bm{h}_1(\bz ,t) 
	=
	& 
	\varrho(z_1) 
	\begin{pmatrix}
		\cos \varphi(z_2) \\ \sin \varphi(z_2) 
	\end{pmatrix} 
	t 
	+ 
	\sum_{n=1}^{s} z_{n+2} \br_n(t),\label{eq:param_2_arcs_1} 
	\\
	\bm{h}_2(\bz ,t) 
	=
	&
	d(z_{s+3})
	\begin{pmatrix}
		\cos \varphi(z_{s+4}) \\ \sin \varphi(z_{s+4}) 
	\end{pmatrix}
	+ 
	\varrho(z_{s+5}) 
	\begin{pmatrix}
		\cos \varphi(z_{s+6}) \\ \sin \varphi(z_{s+6}) 
	\end{pmatrix} 
	t 
	\\
	&
	+ 
	\sum_{n=1}^{s} z_{n+s+6} \br_n(t),\label{eq:param_2_arcs_2}
\end{align}
\end{subequations}
where $d(z) = (d_\text{max}-d_\text{min})(z+\frac{1}{2})+d_\text{min}$, and the functions $\varrho,\varphi$ are the same that in \eqref{eq:geofunsaux}. Finally using the in-variance of the fundamental solution under translations we have that, 
$$
\mathbf{G}(\br^k(\by^k,t),\br^j(\by^j,\tau) = \mathbf{G}(\bm{h}_1(\bz^{k,j},t),\bm{h}_2(\bz^{k,j},\tau)), \quad \forall t,\tau \in [-1,1].
$$
This last equation justify the fact that we only need to construct an interpolation of the function:
$$
\mathbb{H}(x_p,x_q,\bz) = \mathbf{G}(\bm{h}_1(\bz,x_p),\bm{h}_2(\bz,x_q)), \quad p,q =1,\hdots,N_c.
$$
The implementation of the EIM for this particular functions is given in Algorithm \ref{euclid} for the offline part, and Algorithm \ref{eimeval} for the online evaluation\footnote{In the algortims we use the following notation, if $A$ is a two-dimensional array (matrix), its element $(j,k)$, is denoted by $A[j][k]$, its row $j$ is $A[j][:]$, and its column $k$ is $A[:][k]$}. 

\begin{algorithm}
\caption{Empirical Interpolation Method: Offline Phase}\label{euclid}
\begin{algorithmic}[1]
\Procedure{EIM\_OFFLINE\_KERNEL}{$\epsilon_{eim},\bm{z}_s,Q_{\text{max}},QuadPoints$}
\State  $N_s = length(z_s)$
	\For{$\ell=1,\dots,N_s$}
   		\State 
		$H[\ell][:] 
		= 
		\mathbb{H}(QuadPoints,\bm{z}_s(\ell))$ \Comment{Pre-compute the evaluations $\mathbb{H}$}
	\EndFor
	\State $\mathfrak{I} = zeros$ \Comment{Initialize interpolation for all the geometries in $\bm{z_s}$}
	\State $e_{max} = 1e28$ \Comment{The error initialized as a big number}
	\State $\mathfrak{I}_B = []$ \Comment{Stores the interpolation basis }
	\State $\mathfrak{I}_M = [] $
	\Comment{Stores the interaction matrices, $L^R(\mathfrak{I}_B)$}
	\State $X_{max} = []$ \Comment{Stores indices of the quadrature points where big errors are detected}
	\While{$q < Q_{\text{max}}$ and $e_{max} > \epsilon_{eim}$}
	\For{$\ell=1,\dots,N_s$}
		\State 
		Error($\ell$) 
		= 
		$\frac{\|H[\ell][:]-\mathfrak{I}[\ell][:]\|}{\|H[\ell][:]\|}$
	\EndFor
	\State $\ell_\text{max} = \text{argmax}(\text{Error})$   
	\State $\text{e}_\text{max}$ = Error($\ell_\text{max}$) 
	\State 
	$\text{x}_\text{max} =  \text{argmax}(\frac{\|H[\ell_{max}][:]-\mathfrak{I}[\ell_{max}][:]\|}{\|H[\ell_{max}][:]\|})$ 
	\State $X_{max} \gets x_{max}$ 
	\State 
	$\mathfrak{I}_B \gets\frac{H[\ell_{max}][:]-\mathfrak{I}[\ell_{max}][:]}{H[\ell_{max}][x_{max}]-\mathfrak{I}[\ell_{max}][x_{max}]}$	
	\State{$\mathfrak{I}_M \gets L^R\left(\frac{H[\ell_{max}][:]-\mathfrak{I}[\ell_{max}][:]}{H[\ell_{max}][x_{max}]-\mathfrak{I}[\ell_{max}][x_{max}]}\right)$ }
	\State$g = (H[:][X_{max}])^\top$ 
	\State$c = \text{LinSolve}(\mathfrak{I}_B[:][X_{max}],g)$ 
	\State$\mathfrak{I} = (\mathfrak{I}_B)  c$ 
\EndWhile
\State $\mathfrak{I}_s = \mathfrak{I}_B[:][X_{max}]$
\State Return $\mathfrak{I}_S,\mathfrak{I}_M,X_{max}$
\EndProcedure
\end{algorithmic}
\end{algorithm}

\begin{algorithm}
\caption{EIM Evaluation:}\label{eimeval}
\begin{algorithmic}[1]
\Procedure{EIM\_ONLINE}{$QuadPoints,\bm{z},\mathfrak{I}_S,X_{max},\mathfrak{I}_M$}
\State $\bg = \mathbb{H}(QuadPoints[X_{max}],\bm z)$ \Comment{Evaluate the Fundamental solution for the configuration $z$, but only on the specified quadrature points}
\State$\bc = \text{LinSolve}(\mathfrak{I}_S,g)$ \Comment{Obtain the coefficients to exactly interpolate the configuration given by $z$ in the specified Quadrature points}
\State$Interp = \mathfrak{I}_M\bc$ \Comment{Construction of the interpolation}
\State Return $Interp$
\EndProcedure
\end{algorithmic}
\end{algorithm}

\begin{remark}
Lines 19, and 20, of Algorithm \ref{euclid}, ensure that the matrix $\mathfrak{I}_S$, is a triangular matrix, and then corresponding linear system, on Algorithm \ref{eimeval} can be solved fast. Alternatively, we could replace these two lines by their simpler counterparts, 
\begin{align*}
    \mathfrak{I}_B \gets H[\ell_{max}][:]\\
    \mathfrak{I}_M \gets L^R(H[\ell_{max}][:])
\end{align*}
and then instead of returning the square matrix $\mathfrak{I}_B[:][X_\text{max}]$, we return its inverse, so the evaluation in Algorithm \ref{eimeval}, can be carried out fast. This alternative could lead to ill-conditioned linear system in Algorithm \ref{euclid}, however since these are typically very small and solved by direct methods no extra complications arise. 
The main advantage of using this procedure is that in some implementations evaluating \begin{equation}
	L^R
	\left(
		\frac{H[\ell_{max}][:]-\mathfrak{I}[\ell_{max}][:]}{H[\ell_{max}][x_{max}]-\mathfrak{I}[\ell_{max}][x_{max}]}		\right)
\end{equation}
may not be supported, but $ L^R(H[\ell_{max}][:])$ is just a discretization matrix.
\end{remark}

In terms of computational cost, since we are interested in solving for a number of geometric configurations, the only relevant part is Algorithm \ref{eimeval}. For each cross-interaction the cost is $O(Q^2+R^2Q)$. 

Finally, let us comment in the approximation of the self-interaction matrices. As mentioned in Remark \ref{rem:matints}, in an abstract setting the only difference between cross and self-interactions is that for the self-interactions we have two different functions (a regular part, and the function $J(t,\tau)$ of Remark \ref{rem:matints}). 
Thus, we need to interpolate two different functions, but in essence the costs and algorithms are the same.

With this final considerations we have that total cost of the construction of the linear system is $O(M^2(Q^2+R^2Q))$, compared with $O(M^2N^2\log \N)$ for the high-fidelity solver of Section \ref{sec:high_fidelity_solver}.

\section{Analysis of the Reduced Basis Method for Multiple Open Arcs}
\label{sec:convergence_analysis}
In this section, we provide an analysis of the reduced basis method
for the multiple arc problem described in Section \ref{eq:rom_multiple_arcs}.
\subsection{Parametric Holomorphy}
\label{sec:parametric_holomorphy}

We establish the analytic or holomorphic
dependence of the \emph{discrete} parameter-to-solution map
upon the the parametric variables used to describe the arc's shapes. 
This property is of key importance for the derivation of the convergence
analysis presented in Section \ref{sec:rb_convergence}.
The results to be presented herein are based on our previous work \cite{pinto2024shape}.

For $\varrho>1$, we consider the 
Bernstein ellipse in the complex plane
\begin{align}	\mathcal{E}_{\varrho}
	\coloneqq 
	\left\{ 		\frac{\varrho+\varrho^{-1}}{2}:  z\in\IC\;\text{with} \; 1\leq \snorm{z}\leq \varrho
	\right \} 
	\subset \IC.
\end{align}
This ellipse has foci at $z=\pm 1$ and semi-axes of length 
$a\coloneqq  \frac{1}{2}({\varrho}+{\varrho}^{-1})$ and $b \coloneqq  \frac{1}{2}({\varrho}-{\varrho}^{-1})$.
Let us consider the tensorized poly-ellipse
\begin{align}
	\mathcal{E}_{\boldsymbol{\rho}} 
	\coloneqq  
	\bigotimes_{j\geq1} 
	\mathcal{E}_{\rho_j} \subset \IC^{\mathbb{N}},
\end{align}
where $\boldsymbol\rho \coloneqq  \{\rho_j\}_{j\geq1}$
is such that $\rho_j>1$, for $j\in \mathbb{N}$.
\begin{definition}[{\cite[Definition 2.1]{CCS15}}]\label{def:bpe_holomorphy}
Let $X$ be a \emph{complex} Banach space equipped with the norm $\norm{\cdot}_{X}$. 
For $\varepsilon>0$ and $p\in(0,1)$, we say that the map 
\begin{equation}
	\sU \ni \by \mapsto u(\by) \in X
\end{equation}
is $(\boldsymbol{b},p,\varepsilon)$-holomorphic if and only if
\begin{itemize}
	\item[(i)] 
	The map $\sU \ni \by \mapsto u(\by) \in X$ is uniformly bounded, i.e.~
	\begin{align}
		\sup_{\by\in \sU}
		\norm{u(\by)}_{X}
		\leq
		C_0,
	\end{align}
	for some finite constant $C_0>0$.
	\item[(ii)]
	There exists a positive sequence $\boldsymbol{b}\coloneqq \{b_j\}_{j\geq 1} \in \ell^p(\mathbb{N})$ 
	and a constant $C_\varepsilon>0$ such that for any sequence $\boldsymbol\rho\coloneqq \{\rho_j\}_{j\geq1}$ 
	of numbers strictly larger than one that is $(\boldsymbol{b},\varepsilon)$-admissible, i.e.~satisyfing
	\begin{align}
	\label{eq:admissible_polyradius}	
		\sum_{j\geq 1}(\rho_j-1) b_j 
		\leq 
		\varepsilon,
	\end{align}
	the map $\by \mapsto u(\by)$ admits a complex
	extension $\boldsymbol{z} \mapsto u(\boldsymbol{z})$ 
	that is holomorphic with respect to each
	variable $z_j$ on a set of the form 
	\begin{align}
		\mathcal{O}_{\boldsymbol\rho} 
		\coloneqq  
		\displaystyle{\bigotimes_{j\geq 1}} \, \mathcal{O}_{\rho_j},
	\end{align}
	where $\mathcal{O}_{\rho_j}\subset \IC$ is an open set containing $\mathcal{E}_{\rho_j}$.
	This extension is bounded on $\mathcal{E}_{\boldsymbol\rho}$ according to
	\begin{align}
	\label{eq:bpe_hol_bound_epsilon}
		\sup_{\boldsymbol{z}\in \mathcal{E}_{\boldsymbol{\rho}}} 
		\norm{u(\boldsymbol{z})}_{X}  
		\leq 
		C_\varepsilon.
	\end{align}
\end{itemize}
\end{definition}
By considering the parametric  family of open arcs
$\br({\by},\cdot):(-1,1) \rightarrow \IR^2$.
We define the discrete parameter-to-solution map as
$
	\by
	\mapsto
	\boldsymbol{u}_N(\by)
$,
where, for each $\by \in \text{U}$, $\boldsymbol{u}_N(\by) \in \mathbb{T}_N$
corresponds to the solution of \eqref{eq:genbipd1} on the parametrically
defined open arc $\br(\by)$.

\begin{lemma}\label{cor:holomorphy_discrete_d2s}
Let  $m \in \IN$ and $\alpha \in [0,1]$
be such that $m+\alpha>2$. Let Assumpton \eqref{assump:b_decay} be satsfied
with $\boldsymbol{b} \in \ell^p(\IN)$ and $p \in (0,1)$.
Then there exist $\varepsilon>0$ and $N_0 \in \IN$ 
such that for any $N\geq N_0$ the maps
\begin{align}
	\normalfont\text{U}
	\ni
	\by
	\mapsto 
	{\bm u}_N(\by),
	\widetilde{\bm u}_N(\by)
	\in 
	\IT_N
\end{align}
are $(\boldsymbol{b},p,\varepsilon')$-holomorphic
and continuous with the same $\boldsymbol{b} \in \ell^p(\IN)$
and $p\in (0,1)$, and with $\varepsilon'>0$ independent of $N$.
\end{lemma}

\begin{proof}
This result is a direct consequence of stability of the spectral Galerkin BEM discretization,
which based on standard arguments is obtained starting from a base level 
of refinement, thus the validity of ther result for $N\geq N_0$, and
the main result of \cite{pinto2024shape} applied to the BIOs introduced
in Section \ref{sec:model_problem}.
\end{proof}

An equivalent statement for the multiple problem can 
be proved, i.e. parametric holomorphy of the discrete
domain-to-solution map. However, as the ensuing analysis of the ROM
algorithm is based on the understanding of the single arc problem,
we skip it. 

A commonly used concept in nonlinear approximation to 
quantify uniform error bounds is the so-called 
Kolmogorov’s width.
For a compact subset $\mathcal{K}$ of a Banach space
$X$ it is defined for $R\in \IN$ as
\begin{align}
	d_R(\mathcal{K},X)
	\coloneqq
	\inf _{\operatorname{dim}\left(X_R\right) \leq R} 
	\sup _{v \in \mathcal{K}} 
	\min _{w \in X_R}
	\norm{v-w}_X,
\end{align}
where the outer infimum is taken over all
finite dimensional spaces $X_R\subset X$ of  
dimension smaller than $R$. This quantifies
the suitability of $R$-dimensional subspaces
for the approximation of the solution manifold.

In particular, if we consider the response surface
associated to the discrete single arc parametric problem, i.e. $\mathcal{M}^\theta_N$ 
as in 
we have as a consequence of \cite{cohen2016kolmogorov}
and Lemma \ref{cor:holomorphy_discrete_d2s} the following
bound
\begin{equation}
	d_R\left(\mathcal{M}^\theta_N,\IT^{-\frac{1}{2}}\right)
	\leq
	C
	R^{-\left(\frac{1}{p}-1\right)},
	\quad
	R \in \mathbb{N},
\end{equation}
for some $C>0$ and $p \in (0,1)$ as in Assumption \ref{assump:b_decay}.
The exact same result holds valid for the soluton manifold ${\mathcal{M}}^{\theta+\frac{\pi}{2}}_N$.
\subsection{Convergence Analysis}
\label{sec:rb_convergence}

In this section we provide a complete error analysis of the reduced
basis method for multiple open arcs.
The results of this section justify the algorithms from 
Section \ref{eq:rom_multiple_arcs}, which in turn lead
to the numerical results presented ahead in
Section \ref{sec:numerical_results}.

Let $V^{\text{(rb)}}_R$ be as in \eqref{eq:optimal_reduced_space} for $R \in \mathbb{N}$.
Herein, we interested in quantifying the performance
of the RB algorithm according to the following error measure
\begin{equation}\label{eq:error_measure_R}
	\varepsilon\left(V^{\text{(rb)}}_R\right)
	= 
	 \int\limits_{\text{U}}  
	 \hdots  
	 \int\limits_{\text{U}} 
	 \norm{
	 	\mathfrak{u} (\by^1, \hdots, \by^M) 
		- 
		\mathfrak{u}^{\text{(rb)}} _R(\by^1_{\{1:s\}},\hdots \by^M_{\{1:s\}}) 
	}^2_{\IT^{-\frac{1}{2}} \times \cdots \times \IT^{-\frac{1}{2}} } 
	\text{d} \by^1 \cdots \text{d} \by^M,
\end{equation}
where $\by^j_{\{1:s\}} \in \text{U}$, $j=1,\dots,M$, is such that the first $s$ components are equal to the
components of the integraton variable $\by^j$, $j=1,\dots,M$, and the rest of the parametric variables are
set to zero, and where for $\by^1, \hdots, \by^M \in \text{U}$
$$
	\mathfrak{u} (\by^1, \hdots, \by^M)  
	= 
	\begin{pmatrix}
		\bm{u}^1 (\by^1, \hdots, \by^M)  \\
		\vdots \\
		\bm{u}^M (\by^1, \hdots, \by^M)
	\end{pmatrix} 
	\in 
	\IT^{-\frac{1}{2}} \times \cdots \times \IT^{-\frac{1}{2}}
$$
corresponds to the high-fidelity solution of the multiple open arcs problem.
Equivalently, $\mathfrak{u}^{\text{(rb)}} _R\left(\by^1_{\{1:s\}},\hdots \by^M_{\{1:s\}}\right)$ is defined
as the solution of reduced multiple arc problem. 

In the remainder of this section we investigate appropriate bounds for $\varepsilon\left(V^{\text{(rb)}}_R\right)$ in two cases. 
\begin{itemize}
	\item[(i)]
	Firstly, in Section \ref{sec:single_arc_convergence} we
	consider the single arc problem.
	The convergence analysis follows for standard
	arguments, as the ones presented in \cite{hesthaven2016certified,quarteroni2015reduced}.
	A key element in our analysis consists in bounding Kolmogorov's width
	of the solution manifold, which as discussed in \cite{cohen2016kolmogorov}, 
	follows from the parametric holomorphy property of the parameter-to-solution map.
	\item[(ii)]
	Secondly, in Section \ref{sec:multiple_arcs_convergence} we consider the
	multiple open arcs problem.
	The main difficulty is that the reduced basis in this 
	case by construction are only guaranteed to provide good
	performance for the single arc problem. We prove that the multi-arc problem
	can be cast as independent single arc problems but
	on different geometries. 
\end{itemize}

\subsubsection{Convergence Analysis of the Single Arc Problem}
\label{sec:single_arc_convergence}
To simplify the exposition of the multiple arcs case we enumerate
the different steps needed to bound $\varepsilon\left(V^{\text{(rb)}}_R\right)$ 
for the single arc problem with $V^{\text{(rb)}}_R$ as in Section 

\begin{enumerate}
	\item[(i)]  
	{\bf Dimension Truncation in the Parameter Space.}
	As a consequence of the parametric holomorphy
	property of the parameter-to-solution map,
	and following the arguments of \cite{DLC16},
	we truncate the parametric dimension as follows
	\begin{equation}\label{eq:error_quadrature_dim_trunc}
		\varepsilon\left(V^{\text{(rb)}}_R\right)
		\leq
		\int\limits_{\text{U}^{s}}
		\norm{	
			\bm{u}\left(\by_{\{1:s\}}\right)
			-
			\bm{u}^{\text{(rb)}}_R\left(\by_{\{1:s\}}\right)
		}^2_{\IT^{-\frac{1}{2}}} 
		\text{d} \by_s
		+ 
		C_1(p)
		s^{
			-
			2
			\left(
				\frac{1}{p}
				-
				1
			\right)
		}.
	\end{equation}
	for some constant $C_1(p)$ depending on $p\in (0,1)$,
	however not on the parametric dimension $s \in \IN$.
	\item[(ii)]
	{\bf Galerkin BEM Discretization.}
	Recalling the convergence results established in Section \ref{sec:high_fidelity_solver},
	we can replace $\bm{u}\left(\by\right)$ by its 
	Galerkin approximation. 
	Provided that the arc is analyric, as stated in Section \ref{sec:high_fidelity_solver},
	Remark \ref{rmk:convergence_galerkin_bem},
	one has exponential convergence towards the exact solution,
	i.e.~there exist $N_0 \in \IN$, $\varrho>1$ such that for any $N\geq N_0$
	\begin{equation}
	\begin{aligned}
		\varepsilon\left(V^{\text{(rb)}}_R\right)
		\leq 
		&
		\int\limits_{\text{U}^{s}}
		\norm{
			\bm{u}_N\left(\by_{\{1:s\}}\right)
			-
			\bm{u}^{\text{(rb)}}_R\left(\by_{\{1:s\}}\right)
		}^2_{\IT^{-\frac{1}{2}}}
		\text{d} \by_s
		\\
		&
		+
		C_2
		\varrho^{-2N}
		+
		C_1(p)
		s^{
			-
			2
			\left(
				\frac{1}{p}
				-
				1
			\right)
		}.
	\end{aligned}
	\end{equation}
	\item[(iii)]
	{\bf Quasi-optimality in the Reduced Space.}
	Since $\bm{u}^{\text{(rb)}}_R$ is obtained by solving a Galerkin 
	discretization of a well-posed and coercive problem, quasi-optimality 
	yields the existence of a unique discrete solution for $R$
	large enough. I.e., there exists $R_0>0$
	such that for any $R\geq R_0$ one has
	\begin{equation}
	\begin{aligned}
		\varepsilon\left(V^{\text{(rb)}}_R\right)
		\leq
		&
		C_3
		\int\limits_{\text{U}^{(s)}}
		\norm{
			\bm{u}_N
			\left(
				\by_{\{1:s\}}
			\right) 
			-
			P^{\text{(rb)}}_R
			\bm{u}_N
			\left(
				\by_{\{1:s\}}
			\right)
		}^2_{\IT^{-\frac{1}{2}}}   
		\text{d} \by_s
		\\
		&
		+
		C_2
		\varrho^{-2N}
		+
		C_1(p)
		s^{
			-
			2
			\left(
				\frac{1}{p}
				-
				1
			\right)
		},
	\end{aligned}
	\end{equation}
	where $P^{\text{(rb)}}_R: \mathbb{T}_N \rightarrow  V^{\text{(rb)}}_R$
	denotes the orthogonal projection operator
	onto the reduced space $V^{\text{(rb)}}_R$, and $C_3>0$ depends
	only on $R_0$.
	\item[(iv)]
	{\bf Decay of the Kolmmogorov's Width.}
	We are interested in bounding 
	As a consequence of Lemma \ref{cor:holomorphy_discrete_d2s}, which in turn
	follows from our previous work \cite{pinto2024shape} and the stability of the Galerkin 
	BEM discetrization described in Section \ref{sec:high_fidelity_solver}, and \cite{cohen2016kolmogorov}
	we may conclude that for $R\geq R_0$ and $N\geq N_0$, with $R_0$ as in item (iii) and 
	$N_0$ as in item (ii), we have
	\begin{equation}
		\varepsilon\left(V^{\text{(rb)}}_R\right)
		\leq
		C_4
		R^{
			-
			2
			\left(
				\frac{1}{p}
				-
				1
			\right)
		}
		+
		C_2
		\varrho^{-2N}
		+
		C_1(p)
		s^{
			-
			2
			\left(
				\frac{1}{p}
				-
				1
			\right)
		},
	\end{equation}
	for some $C_4>0$.
\end{enumerate}

Clearly, the computation of $V^{\text{(rb)}}_R$ is not feasible as it requires
the evaluation of \eqref{eq:optimal_reduced_space}. Therefore, we consider
a discrete approximation of the high-dimensional integral as in \eqref{eq:minimization_rom_discrete},
and we refer to this space as $\widetilde{V}^{\text{(rb)}}_R$.
To quantify this missmatch, we follow \cite[Section 6.5]{quarteroni2015reduced}. 

Set
\begin{equation}
	\varepsilon^{(N_s)}
	\left(
		\widetilde{V}^{\text{(rb)}}_R
	\right)
	=
	\frac{1}{N_t}
	\sum_{j=1}^{N_t}
	\norm{
		\bm{u}_N
		\left(
			\by_j
		\right) 
		- 
		P^{\text{(rb)}}_R
		\bm{u}_N
		\left(
			\by_j
		\right)
	}^2_{\IT^{-\frac{1}{2}}}.
\end{equation}
We are interested in the performance of this empirically constructed space 
according to
\begin{equation}
	\varepsilon
	\left(
		\widetilde{V}^{\text{(rb)}}_R
	\right)
	\leq
	\underbrace{
	\snorm{
		\varepsilon
		\left(
			\widetilde{V}^{\text{(rb)}}_R
		\right)
		-
		\varepsilon^{(N_s)}
		\left(
			\widetilde{V}^{\text{(rb)}}_R
		\right)
	}}_{(\clubsuit)}
	+
	\varepsilon^{(N_s)}
	\left(
		\widetilde{V}^{\text{(rb)}}_R
	\right)
\end{equation}
It follows from \cite{halton1960efficiency} that for each $\delta>0$ there exists $C(\delta)>0$
such that $(\clubsuit) \leq C(\delta) N^{-1+\delta}$, and with $C(\delta) \rightarrow \infty$
as $\delta\rightarrow 0^+$.
The last term can be bounded by the singular values of the snapshot matrix.

\subsubsection{Convergence Analysis for Multiple Open Arcs}
\label{sec:multiple_arcs_convergence}
We proceed to establish the convergence of the RB method 
applied to the multiple arcs problem as described in
Section \ref{sec:reduced_basis_multiple_arc}.
Throughout, we work under the following assumption.

\begin{assumption}\label{assump:rotations_arcs}
The set $\Sigma$ defined in \eqref{eq:collection_of_arcs} is closed
under rotations, i.e.~if $\br \in \Sigma$,
then for any $\theta \in [0,2\pi)$, $\bm{R}_\theta \br \in \Sigma$,
where $\bm{R}_\theta$ denotes the rotation matrix for angle $\theta$
defined in \eqref{eq:rotmatrix}.
In addition, for $s \in \IN$, $s> 4$, the set $\Sigma_{s}$ is
closed under rotations as well.
\end{assumption}

\begin{remark}
\hfill
\begin{itemize}
	\item[(i)]
	For $s = 4$ is immediate that $\Sigma_{s}$ is closed
	under rotations as they are only line segments in $\IR^2$.
	\item[(ii)]
	For larger values of $s$, i.e. $s>4$, in general it would depend 
	on the properties of the functions $\{\br_n\}_{n\in \IN}$.
	\item[(iii)]
	A particular case for which one may straightforwardly verify 
	Assumption \ref{assump:rotations_arcs} is when the
	functions $\{\br_n\}_{n\in \IN}$ are of the form $r_n \bm{e}_1$, $r_n \bm{e}_2$,
	for some scalar functions $r_n$, and for each $n \in \IN$.
\end{itemize}
\end{remark}

For the sake of simplicity, we firstly discuss the two open arcs
problem as the extension to multiple open arcs follows from the exact
same arguments.
This problem reads as follows:
For $\by^1,\by^2 \in \text{U}$, we seek
$\bm{u}^1(\by^1,\by^2) , \bm{u}^2(\by^1,\by^2)  \in \IT^{-\frac{1}{2}}$
such that
\begin{equation}
\begin{aligned}
	\bm{V}_{\br^1(\by ^1),\br^1(\by^1)} \bm{u}^1(\by^1,\by^2) 
	+ 
	\bm{V}_{\br^1(\by^1),\br^2(\by^2)}\bm{u}^2(\by^1,\by^2) & = \bm{g}_1,\\
        \bm{V}_{\br^2(\by^2),\br^1(\by^1)}\bm{u}^1(\by^1,\by^2)
        +         \bm{V}_{\br^2(\by^2),\br^2(\by^2)} \bm{u}^2(\by^1,\by^2)  & = \bm{g}_2,
\end{aligned}
\end{equation} 
where $\bm{g}_1,\bm{g}_2$ are as in \eqref{eq:rhsnot}.

The analysis of the multiple arc problem
has one major difference compered to the single arc case:
In the former case the  cross-interaction terms cannot
be acounted by the use of the reduced space.


In the following we argue that the solution of the
multiple open arc problem can be approximated
(in a controlled way) by a linear combinations of
independent problem, posed on different geometries.
First notice that the multi-arc problem can be expresed as 
\begin{align}
	\bm{V}_{\br^1(\by^1),\br^1(\by^1)} \bm{u}^1\left(\by^1,\by^2\right) = \bm{f}_1  \left(\by^1,\by^2 \right)
	\quad
	\text{and}
	\quad	\bm{V}_{\br^2(\by^2),\br^2(\by^2)} \bm{u}^2 \left(\by^1,\by^2 \right) = \bm{f}_2  \left(\by^1,\by^2 \right),
\end{align}
where 
\begin{equation}
\begin{aligned}
	\bm{f}_1 \left(\by^1,\by^2 \right)  &  = \bm{g}_1 - \bm{V}_{\br^1(\by^1),\br^2(\by^2)} \bm{u}^2 \left(\by^1,\by^2 \right), \;\; \text{and},  \\
	\bm{f}_2 \left(\by^1,\by^2 \right)  &  = \bm{g}_2 - \bm{V}_{\br^2(\by^2),\br^1(\by^1)} \bm{u}^1\left(\by^1,\by^2 \right) 
\end{aligned}
\end{equation} 

Next, we approximate $\bm{f}_1,$ and $\bm{f}_2$.
To this end, we recall that the set of boundary traces of plane-waves (of a fixed wave-number and variable directions) are dense in $L^2(\gamma)$, where $\gamma$ is the boundary of a bounded, simply connected Lipschitz domain, provided that the wavenumber is not a eigenvalue
of the interior Laplace Dirichlet problem, see e.g.~\cite[Section 3.4]{colton1998inverse}.
Standard arguments for open arcs (see, e.g., \cite{Sauter:2011})
guarantee that we can extend this result to the scenario of $\gamma$ being an open arc.

%

Being $L^2(\gamma)$ a dense subset of $\mathbb{W}^{\frac{1}{2}}(\gamma)$, 
one may claim the following: 
For each $\by^1,\by^2 \in \text{U}$ and $\epsilon >0$ there exist $L \in \IN$ and 
\begin{equation}\label{eq:alpha_beta_theta}
	\quad
	\left\{\alpha_1^\ell(\by^1,\by^2) \right\}_{\ell=1}^{L}\subset \IC, 
	\quad
	\left\{\beta_1^\ell(\by^1,\by^2) \right\}_{\ell=1}^{L}\subset \IC,
	\quad
	\text{and}
	\quad
	\left\{\theta^\ell_1(\by^1,\by^2)\right\}_{\ell=1}^{L} \subset [0,2\pi),
\end{equation}
such that
\begin{equation}
	\norm{
		\bm{f}_1  
		\left(\by^1,\by^2 \right)
		- 
		\sum_{\ell=1}^L 
		\alpha^\ell_1
		\bm{e}_{\theta^\ell_1}g_{\theta^\ell_1,k_\mathsf{p}} \circ \br^1(\by^1)
		+ 
		\beta^\ell_1 
		\bm{e}_{
			\theta^\ell_1
			+ 
			\frac{\pi}{2}
		}
		g_{\theta^\ell_1,k_\mathsf{p}} \circ \br^1(\by^1) 
	}_{\IW^{\frac{1}{2}}}
	< 
	\epsilon,
\end{equation}
and a equivalent result holds for $\bm{f}_2\left(\by_1,\by_2 \right)$.
We remark that the quantities in \eqref{eq:alpha_beta_theta}
depend continuously on the parametric inputs $\by^1,\by^2 \in \text{U}$ and on $\epsilon>0$.
For $\by^1, \by^2 \in \text{U}$, let us now define the collection of functions
$$
	\left\{\bm{v}^\ell_j(\by^1,\by^2) \right\}_{\ell=1}^{L}
	\subset \IT^{-\frac{1}{2}},
	\quad
	\text{and}
	\quad
	\left\{\bm{w}^\ell_j(\by^1,\by^2)  \right\}_{\ell=1}^{L} \subset \IT^{-\frac{1}{2}},
$$
for $j=1,2$, as
\begin{subequations}
\begin{align}
	\bm{V}_{\br^j(\by^j),\br^j(\by^j)} \bm{v}^\ell_j(\by^1,\by^2) 
	&
	= 
	\bm{e}_{\theta_j^l} g_{\theta_j^l,\kappa_\mathsf{p}} \circ  \br^j(\by^j)
	\quad
	\text{and}, \\
	\bm{V}_{\br^j(\by^j),\br^j(\by^j)} \bm{w}^\ell_j (\by^1,\by^2) 
	&
	= 
	\bm{e}_{\theta_j^l +\frac{\pi}{2}} 
	g_{\theta_j^l,\kappa_\mathsf{p}} \circ  \br^j(\by^j),
\end{align}
\end{subequations}
and their rotations by the angle $\phi^\ell_j (\by^1,\by^2) = \theta - \theta_j^\ell(\by^1,\by^2) $,
denoted by
\begin{equation}
	\widetilde{\bm{v}}^\ell_j(\by^1,\by^2) = \bm{R}_{\phi^\ell_j(\by^1,\by^2) } \bm{v}^\ell_j(\by^1,\by^2)
	\quad
	\text{and}
	\quad
	\widetilde{\bm{w}}^\ell_j(\by^1,\by^2) = \bm{R}_{\phi^\ell_j(\by^1,\by^2) } \bm{w}^\ell_j(\by^1,\by^2),
\end{equation}
for $j=1,2$.

Given $\theta, \tilde{\theta} \in [0,2\pi)$ one has 
$g_{\theta,\kappa}(\bm{R}_{\theta-\tilde{\theta}} \bx ) = g_{\tilde{\theta},\kappa}(\bx)$.
Also, observe that $\bm{G}$ is invariant under translations,~i.e.
$\bm{G}(\bx+ \bm{d} , \by +\bm{d}) = \bm{G}(\bx,\by)$ for any
$\bx,\by,\bm{d}\in \mathbb{R}^2$.
However, it is not invariant under rotations. Instead, the following holds
$
	\bm{G}(\bm{R}_\theta \bm{x},\bm{R}_\theta \bm{y})
	= 
	\bm{R}_\theta \bm{G}(\bm{x},\bm{y}) \bm{R}_\theta^\top.
$
These observations imply that
\begin{equation}
\begin{aligned}
	\bm{V}_{\widetilde{\br}^j(\by^j),\widetilde{\br}^j(\by^j)} 
	\widetilde{\bm{v}}^\ell_j 
	(\by^1,\by^2) 
	&
	= 
	\bm{e}_{\theta} g_{\theta,\kappa_\mathsf{p}} 
	\circ  
	\widetilde{\br}^j(\by^j)
	\quad
	\text{and}
	\\
	\bm{V}_{\widetilde{\br}^j(\by^j),\widetilde{\br}^j(\by^j)} 
	\widetilde{\bm{w}}^\ell_j 
	(\by^1,\by^2) 
	&
	=  
	\bm{e}_{\theta+\frac{\pi}{2}} g_{\theta,\kappa_\mathsf{p}} 
	\circ  
	\widetilde{\br}^j(\by^j),
\end{aligned}
\end{equation} 
where $\widetilde{\br}^j(\by^j) = \bm{R}_{\phi^\ell_j} \br^j(\by^j)$, for $j = 1,2$.

It follows from Assumption \ref{assump:rotations_arcs} we have that  $\widetilde{\br}^j(\by^j)$
belongs to $\Sigma$. 
Consequently, $ \widetilde{\bm{v}}^\ell_j,\widetilde{\bm{w}}^\ell_j$
are in the solution manifold of the single arc problem.
In this setting, we obtain the following approximation result. 

\begin{proposition}\label{prop:approxsols}
Let $V$ be any subspace of $\IT^{-\frac{1}{2}}$.
Given $\epsilon >0$ there exist $L \in \IN$, $\{\theta^\ell_j\}_{\ell=1}^L \subset [0,2\pi)$
and functions $\{\alpha^\ell_j(\by^1,\by^2)\}_{\ell=1}^L,
\{\beta_j^\ell(\by^1,\by^2)\}_{\ell=1}^L \subset \IC$, depending on the parametrc inputs
$\by^1,\by^2$, such that for $j=1,2$
\begin{equation}\label{eq:error_rotation_total}
\begin{aligned}
	\norm{
		\bm{u}^j
		(\by^j)
		- 
		P_{V}
		\bm{u}^j
		(\by^j)
	}_{\IT^{-\frac{1}{2}}} 
	\\
	&
	\hspace{-2cm} 
	\lesssim
	\epsilon 
	+ 
	 \sum_{\ell=1}^L 
	\snorm{\alpha_j^\ell(\by^1,\by^2)}
	\left(\norm{
		\widetilde{\bm{v}}_j^\ell(\by^1,\by^2) 
		- 
		P_{V}
		\widetilde{\bm{v}}_j^\ell(\by^1,\by^2) 
	}_{\IT^{-\frac{1}{2}}} 
	+
	\mathcal{E}
	\left(
		\mathcal{M}, V
	\right)
	\right)
	\\
	&
	\hspace{-1.7cm} 
	+ 
	\sum_{\ell=1}^L 
	\snorm{\beta_j^\ell(\by^1,\by^2)}
	 \left(
	\norm{
		\widetilde{\bm{w}}_j^\ell(\by^1,\by^2) 
		- 
		P_{V}		\widetilde{\bm{w}}_j^\ell(\by^1,\by^2) 
	}_{\IT^{-\frac{1}{2}}}
 	+
	\mathcal{E}
	\left(
		\widetilde{\mathcal{M}}, V
	\right)
 \right),
\end{aligned}
\end{equation}
where 
\begin{equation}\label{eq:rotation_error}
	\mathcal{E}
	\left(
		\mathcal{M},V
	\right)
	= 
	\sup_{\substack{\theta \in [0,2\pi) \\ \bm{x} \in P_V( \mathcal{M})}} 
	\norm{
		\bm{R}_\theta \bm{x} 
		- 
		P_{V}
		\bm{R}_\theta 
		\bm{x} 
	}_{\IT^{-\frac{1}{2}}}
\end{equation}
measures how well rotations of the set $\mathcal{M}$
can be approximated by the space $V$.
\end{proposition}

\begin{proof}
First, for $j=1,2$ and each $\by^j \in \text{U}$ observe that 
\begin{equation}\label{eq:triangle_ineq_error}
\begin{aligned}
	\norm{
		\bm{u}^j 
		(\by^j)
		-
		P_{V}
		\bm{u}^j
		(\by^j)
	}_{\IT^{-\frac{1}{2}}}
	\leq 
	&
	2\norm{
		\bm{u}^j(\by^j)
		- 
		\sum_{\ell=1}^L 
		\alpha_j^\ell(\by^1,\by^2) \bm{v}_j^\ell
		+ 
		\beta_j^\ell(\by^1,\by^2) \bm{w}_j^\ell
	}_{\IT^{-\frac{1}{2}}}
	\\
	&
	+ 
	\sum_{\ell=1}^L 
	\snorm{
		\alpha_j^\ell
		(\by^1,\by^2)
	}
	\norm{
		\bm{v}_j^\ell(\by^1,\by^2) 
		- 
		P_{V}
		\bm{v}_j^\ell(\by^1,\by^2) 
	}_{\IT^{-\frac{1}{2}}} 
	\\
	&
	+ 
	\sum_{\ell=1}^L 
	\snorm{\beta_j^\ell(\by^1,\by^2)}
	\norm{
		\bm{w}_j^\ell (\by^1,\by^2) 
		- 
		P_{V}
		\bm{w}_j^\ell(\by^1,\by^2) 
	}_{\IT^{-\frac{1}{2}}}.
\end{aligned}
\end{equation}
In the following, for economy of notation, we dot explictely state the dependence
upon the parametrc variables $\by^1$ and $\by^2$.
The first term on the right-hand side 
of \eqref{eq:triangle_ineq_error} 
can be bounded using the well-possedness
of the boundary integral formulation.
Indeed, for $j=1,2$, one has
\begin{equation}
\begin{aligned}
	\norm{
		\bm{u}^j
		- 
		\sum_{\ell=1}^L 
		\alpha_j^\ell \bm{v}_j^\ell  + \beta_j^\ell\bm{w}_j^\ell
	}_{\IT^{-\frac{1}{2}}} 
	&
	\lesssim
	\norm{
		\bm{V}_{\br^j,\br^j}
		\left( 
 			\bm{u}^j
			- 
			\sum_{l=1}^L \alpha_j^\ell \bm{v}_j^\ell  + \beta_j^\ell \bm{w}^\ell_j
		\right)
	}_{\IW^{\frac{1}{2}}} 
	\\
	&
	= 
	\norm{
		\bm{f}_j
		- 
		\sum_{l=1}^L \alpha^\ell_j \bm{e}_{\theta^\ell_j}g_{\theta^\ell_j,k_p} \circ \br^j
		+ 
		\beta^\ell_j\bm{e}_{\theta^\ell_j+ \frac{\pi}{2}}g_{\theta^\ell_j,k_p} \circ \br^j 
	}_{\IW^{\frac{1}{2}}}
	\\
	&
	<
	\epsilon.
\end{aligned}
\end{equation}

We bound the term
$\norm{\bm{v}_j^\ell - P_{V} \bm{v}_j^\ell}_{\IT^{-\frac{1}{2}}}$, $\ell = 1,\dots, L$, $j=1,2$,
as all the remainder terms are similar. 
Using the definition of $\widetilde{\boldsymbol{v}}_j^l$ and
the invariance of the $2$-norm under rotations we have that
\begin{equation}
\begin{aligned}
	\norm{
		\bm{v}_j^\ell
		- 
		P_{V} \bm{v}_j^\ell
	}_{\IT^{-\frac{1}{2}}} 
	&
	= 
	\norm{
		\bm{R}^\top_{\phi_j^\ell}
		\widetilde{\boldsymbol{v}}_j^\ell 
		- 
		P_{V}
		\bm{R}^\top_{\phi_j^\ell}
		\widetilde{\boldsymbol{v}}_j^\ell
	}_{\IT^{-\frac{1}{2}}} 
	\\
	& 
	\leq 
	\norm{
		\widetilde{\boldsymbol{v}}_j^\ell 
		-
		P_{V} 
		\widetilde{\boldsymbol{v}}_j^\ell
	}_{\IT^{-\frac{1}{2}}} 
	+ 
	\norm{
		\left(
			\bm{R}^\top_{\phi_j^\ell}
			P_{V}
			- 
			P_{V}
			\bm{R}^\top_{\phi_j^\ell}
		\right)
		\widetilde{\boldsymbol{v}}_j^\ell
	}_{\IT^{-\frac{1}{2}}}.
\end{aligned}
\end{equation}
In addition, we have
\begin{equation}
\begin{aligned}
	\norm{
		\left(
			\bm{R}^\top_{\phi_j^\ell}
			P_{V}
			- 
			P_{V}
			\bm{R}^\top_{\phi_j^\ell}
		\right)
		\widetilde{\boldsymbol{v}}_j^\ell
	}_{\IT^{-\frac{1}{2}}}
	\leq
	&
	\norm{
		\left(
			\bm{R}^\top_{\phi_j^\ell}
			P_{V} 
			- 
			P_{V}
			\bm{R}^\top_{\phi_j^\ell}
			P_{V} 
		\right) 
		\widetilde{\boldsymbol{v}}_j^\ell
	}_{\IT^{-\frac{1}{2}}}
	\\
	& 
	+ 
	\norm{
		\left(
			P_{V}  \bm{R}^\top_{\phi_j^\ell} P_{V}
			-
			P_{V} \bm{R}^\top_{\phi_j^\ell}
		\right)
		\widetilde{\boldsymbol{v}}_j^\ell
	}_{\IT^{-\frac{1}{2}}}.
\end{aligned}
\end{equation}
Next, set $\bm{h} = P_{V} \widetilde{\boldsymbol{v}}_j^\ell$.
Then, one has
\begin{equation}
\begin{aligned}
	\norm{
		\left(
			\bm{R}^\top_{\phi_j^\ell}
			P_{V}
			- 
			P_{V}
			\bm{R}^\top_{\phi_j^\ell}
			P_{V}
		\right) 
		\widetilde{\boldsymbol{v}}_j^\ell
	}_{\IT^{-\frac{1}{2}}} 
	&
	= 
	\norm{
		\left(
			\bm{R}^\top_{\phi_j^\ell}
			-	 
			P_{V}
			\bm{R}^\top_{\phi_j^\ell}
		\right) 
		\bm{h}
	}_{\IT^{-\frac{1}{2}}}  
	\\
	&
	\leq 
	\mathcal{E}(\mathcal{M}, V),
\end{aligned}
\end{equation}
where $\varepsilon(\mathcal{M}, V)$ is as in \eqref{eq:rotation_error}.

Using the continuity of $P_{V} \bm{R}^\top_{\phi_j^\ell}$ we obtain
\begin{equation}
	\norm{
		\left(
			P_{V}  \bm{R}^\top_{\phi_j^\ell} P_{V}
			-
			P_{V} \bm{R}^\top_{\phi_j^\ell}
		\right)
		\widetilde{\boldsymbol{v}}_j^\ell
	}_{\IT^{-\frac{1}{2}}}
	\leq 
	\norm{
		P_{V}\widetilde{\boldsymbol{v}}_j^\ell 
		-
		\widetilde{\boldsymbol{v}}_j^\ell
	}_{\IT^{-\frac{1}{2}}},
\end{equation}
thus yielding
\begin{equation}
	\norm{
		\left(
			\bm{R}^\top_{\phi_j^\ell}
			P_{V} 
			- 
			P_{V} 
			\bm{R}^\top_{\phi_j^\ell}
			P_{V} 
		\right) 
		\widetilde{\boldsymbol{v}}_j^\ell
	}_{\IT^{-\frac{1}{2}}}
	\leq 
	\mathcal{E}(\mathcal{M}, \IT^R_N) 
	+ 
	\norm{
		P_{V}
		\widetilde{\boldsymbol{v}}_j^\ell
		-
		\widetilde{\boldsymbol{v}}_j^l
	}_{\IT^{-\frac{1}{2}}},
\end{equation}
and finally replacing in the original bound of
$\norm{\bm{v}^\ell_j - P_V \bm{v}^\ell_j}_{\IT^{-\frac{1}{2}}}$ we obtain the bound
\begin{equation}
	\norm{
		\bm{v}^\ell_j 
		- 
		P_{V} \bm{v}^\ell_j
	}_{\IT^{-\frac{1}{2}}}
	\lesssim 
	\norm{
		\widetilde{\boldsymbol{v}}_j^\ell
		-
		P_{V}
		\widetilde{\boldsymbol{v}}_j^\ell
	}_{\IT^{-\frac{1}{2}}}
	+
	\mathcal{E}
	\left(
		\mathcal{M}
		, 
		V
	\right).
\end{equation}
This allow us to state the bound in \eqref{eq:error_rotation_total}, thus concluding the proof.
\end{proof}

Equipped with this we may now state the main result of this section concerning
the convergence of the RB algorithm from for the multiple arcs problem.
\begin{theorem}\label{teo:convgmult}
Let Assumptions \ref{assump:b_decay} and \ref{assump:rotations_arcs} be satisfied for some $p \in (0,1)$.
There exist $R_0, L,N_0 \in \IN$ and $\varrho>1$ such that for $R\geq R_0$, $N\geq N_0$ it holds
\begin{align*}
	\varepsilon\left(V^{\normalfont\text{(rb)}}_R\right)
	\lesssim	
	\epsilon +s^{-2\left(\frac{1}{p}-1\right)}+ \varrho^{-2N}
	+
	R^{-2\left(\frac{1}{p}-1\right)}
	+ 
	\mathcal{E}^2
	\left(
		\mathcal{M}^\theta_N, V^{\normalfont\text{(rb)}}_R
	\right)	
	+ 
	\mathcal{E}^2
	\left(
		\mathcal{M}^{\theta+\frac{\pi}{2}}_N, V^{\normalfont\text{(rb)}}_R
	\right)
\end{align*}
where $\mathcal{E}\left(\mathcal{M}^\theta_N,V^{\normalfont\text{(rb)}}_R\right)$
and $\mathcal{E}\left(\mathcal{M}^{\theta+\frac{\pi}{2}}_N,V^{\normalfont\text{(rb)}}_R\right)$
are as in Proposition \ref{prop:approxsols},
and $\varepsilon\left(V^{\normalfont\text{(rb)}}_R\right)$ as in \eqref{eq:error_measure_R}.
\end{theorem}

\begin{proof}
By using the same arguments of the analysis for a single arc we obtain the following bound 
{\small
\begin{align*}   
	&
	\varepsilon\left(V^{\normalfont\text{(rb)}}_R\right) 
	\lesssim 
	M \left( s^{-2\left(\frac{1}{p}-1\right)}+ \varrho^{-2N}\right)  
	\\ 
	& 
	+
	\sum_{j=1}^M
	\int_{\text{U}^{s} }\cdots \int_{\text{U}^{s}}  
	\norm{
		\bm{u}_N^j(\by^1_{\{1:s\}},\hdots,\by^M_{\{1:s\}}) 
		- 
		P^{\normalfont\text{(rb)}}_R
		\bm{u}_N^j(\by^1_{\{1:s\}},\hdots,\by^M_{\{1:s\}}) 
	}_{\IT^{-\frac{1}{2}}}^2 
	\text{d} \by^1_{\{1:s\}} \hdots \text{d} \by^M_{\{1:s\}}
\end{align*}
}%
An inspection of Proposition \ref{prop:approxsols} reveals that the stament continues to be valid
for the discrete counterpart of the multiple arc problem.
Consequently, for any $\epsilon>0$ there exist $L \in \IN$ such that for $j = 1,\dots, M$ it holds
\begin{equation}\label{eq:error_bound_M_arcs}
\begin{aligned}    
	\norm{
		\bm{u}_N^j(\by^1_{\{1:s\}},\hdots,\by^M_{\{1:s\}}) - P^{\normalfont\text{(rb)}}_R\bm{u}_N^j(\by^1_{\{1:s\}},\hdots,\by^M_{\{1:s\}})
	}_{\IT^{-\frac{1}{2}}} 
	\\
	&
	\hspace{-6cm}
	\lesssim
	\epsilon 
	+
	\sum_{\ell=1}^L 
	\snorm{\alpha_j^\ell}
	\left(\norm{
		\widetilde{\bm{v}}_j^\ell
		- 
		P^{\normalfont\text{(rb)}}_R
		\widetilde{\bm{v}}_j^\ell
	}_{\IT^{-\frac{1}{2}}} 
	+
	\mathcal{E}
	\left(
		\mathcal{M}^\theta_N, V^{\normalfont\text{(rb)}}_R
	\right)
	 \right)
	&
	\\
	&
	\hspace{-5.7cm}
	+ 
	\sum_{\ell=1}^L 
	\snorm{\beta_j^\ell}
	 \left(
	\norm{
		\widetilde{\bm{w}}_j^\ell
		- 
		P^{\normalfont\text{(rb)}}_R
		\widetilde{\bm{w}}_j^\ell
	}_{\IT^{-\frac{1}{2}}}
 	+
	\mathcal{E}
	\left(
		\mathcal{M}^{\theta+\frac{\pi}{2}}_N, V^{\normalfont\text{(rb)}}_R
	\right)
 \right).
\end{aligned}
\end{equation}
We remark that in \eqref{eq:error_bound_M_arcs}, just for economy of notation,
we did not explictely include the dependence upon the parametric variables $\by^1_{\{1:s\}},\hdots,\by^M_{\{1:s\}}$
in the last two terms.

By Assumption \ref{assump:rotations_arcs} we have that $\widetilde{\bm{v}}^\ell_j$, and $\widetilde{\bm{w}}^\ell_j$
are elements of the discrete solution manifold.
Arguing as in the case of a single arc, in particular as in items (iii) and (iv) of the results presented in
Secton \ref{sec:single_arc_convergence}, we have the following bounds
\begin{align*}
    \norm{
		\widetilde{\bm{v}}_j^\ell
		\left(\by^1_{\{1:s\}},\hdots,\by^M_{\{1:s\}}\right)
		- 
		P^{\normalfont\text{(rb)}}_R
		\widetilde{\bm{v}}_j^\ell
		\left(\by^1_{\{1:s\}},\hdots,\by^M_{\{1:s\}}\right)
	}_{\IT^{-\frac{1}{2}}} & \lesssim R^{- \left(\frac{1}{p}-1\right)}
	\quad
	\text{and}
	\\
   	\norm{
		\widetilde{\bm{w}}_j^\ell
		\left(\by^1_{\{1:s\}},\hdots,\by^M_{\{1:s\}}\right)
		- 
		P^{\normalfont\text{(rb)}}_R
		\widetilde{\bm{w}}_j^\ell
		\left(\by^1_{\{1:s\}},\hdots,\by^M_{\{1:s\}}\right)
	}_{\IT^{-\frac{1}{2}}}  & \lesssim R^{- \left(\frac{1}{p}-1\right)}.
\end{align*}
Finally, it follows from a compactness argument and the continuity of the functions $\alpha^l_j$ and $\beta^l_j$ for a given $\epsilon>0$ the corresponding value of $L$ can be selected such that do not depend of the parameters $\by^1_{\{1:s\}},\hdots,\by^M_{\{1:s\}}$, thus
yielding the stated result.
\end{proof}

Various remarks are in place regarding the convergence of the multiple arcs problem. 
First, the presence of the term $\mathcal{E}\left(\mathcal{M}^\theta_N, V^{\normalfont\text{(rb)}}_R\right) $ in Proposition \ref{prop:approxsols}
is a result of the lack of invariance under rotation of the fundamental solution of the elastic wave problem.
For acoustic formulation or EFIE formulation for Maxwell equation this term does not appear.   
In addition, the value of $L$ and $\epsilon$ of Theorem \ref{teo:convgmult} are obviously related, as $\epsilon>0$ can be interpreted as the convergence rate of the approximation by plane waves of a particular solution of the Navier equation. In particular for the Helmholtz equation the relation is well understood, see e.g.~\cite{moiola2011plane}.

\section{Numerical Results}
\label{sec:numerical_results}
We present numerical experiments illustrating the performance of 
the reduced order algorithm for the multiple arcs problem.

\subsection{Fixed Number of Arcs}
\begin{figure}
    \centering    \includegraphics[scale=0.18]{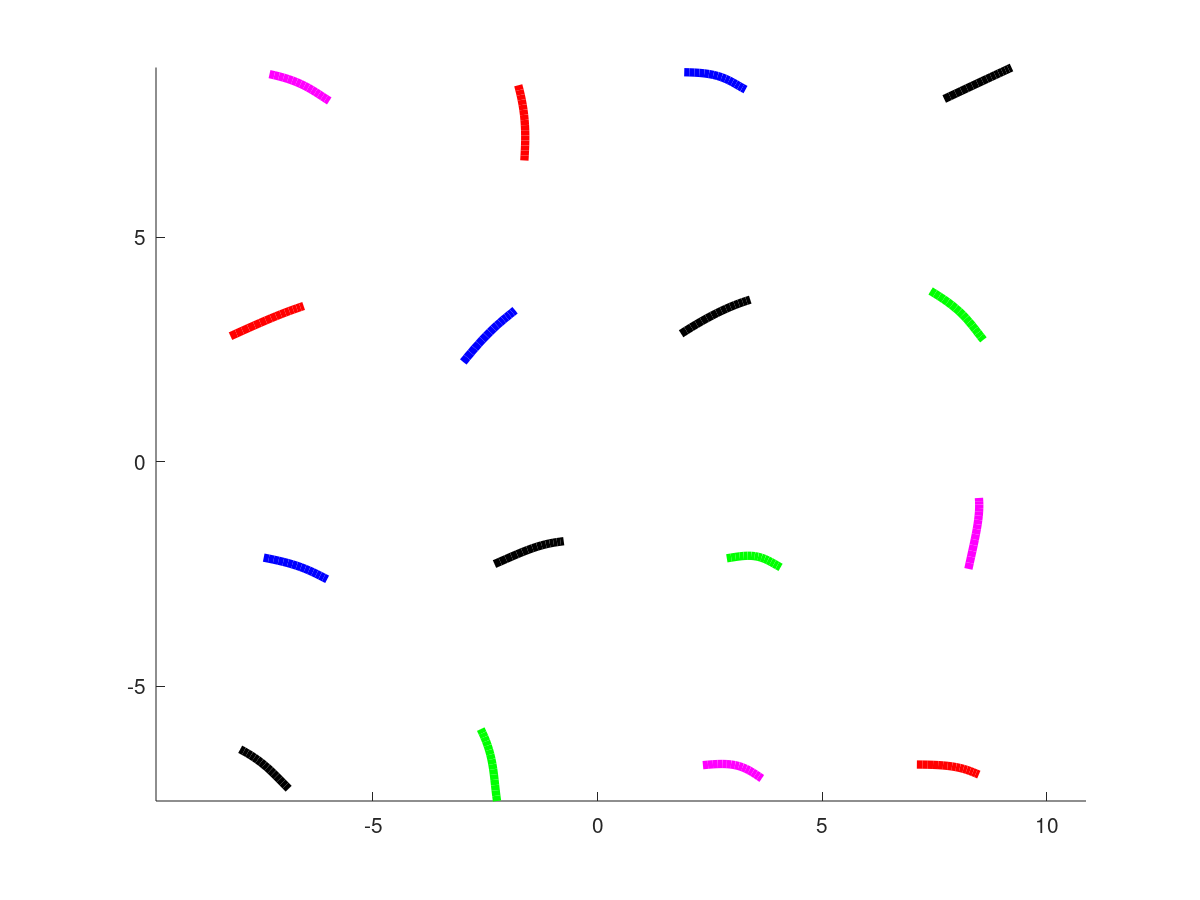}
    \caption{Geometry realization of $16$ open arcs.}
    \label{fig:geo16}
\end{figure}

We consider 16 open arcs enclosed in the two-dimensional box
$[-10,10]\times [-10,10]$. In addition, we consider the setting described
in Section \ref{sec:geo} with parameters
$r_{\text{min}} = 0.56$, $r_{\text{max}} = 0.93$, $d_{\text{min}} = 5$, 
$d_{\text{max}} = 21$, and perturbations terms in 
\eqref{eq:param_2_arcs_1}--\eqref{eq:param_2_arcs_2} of the form
\begin{equation}\label{eq:pert16}
	\sum_{n=1}^3 
	c_n
	\sum_{p=1}^2 
	\left(
		\cos((n-1)t) \bm{e}_p  y_{6(p-1)+2n} 
		+ 
		\sin(nt) \bm{e}_p y_{6(p-1)+2n-1}
	\right), 
	\quad 
	\by \in \text{U}^s,  
\end{equation}
with dimension truncation $s=12$
and $c_n = n^{-\frac{5}{2}}$.
A realization of this setting is presented in Figure \ref{fig:geo16}. 
We consider elastic-wave scattering operator with parameters
$\omega = 10,$ $\lambda =2$, $\mu = 1$.

\begin{figure}
    \centering    \includegraphics[scale=0.18]{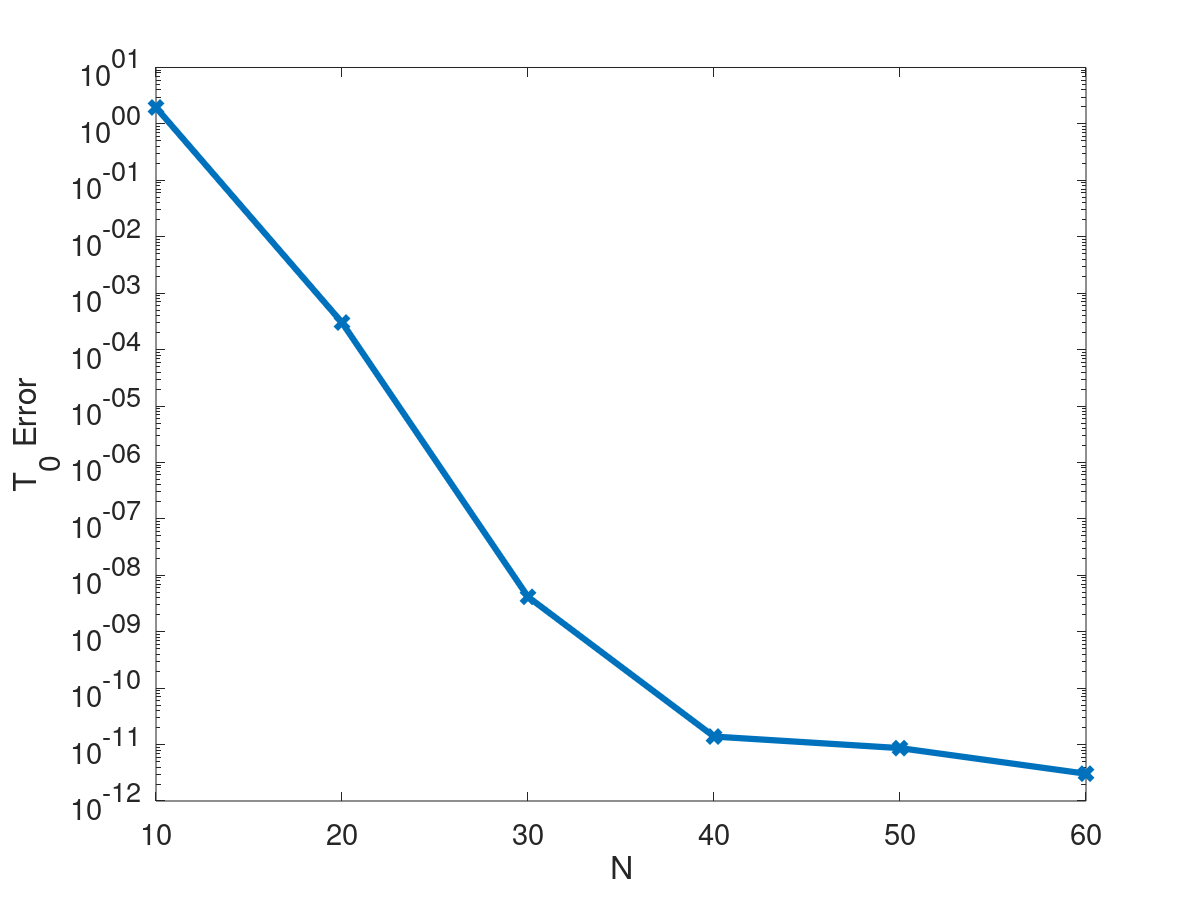}
    \caption{Convergence of the high-fidelity solver.}
    \label{fig:hferrors16}
\end{figure}

We investigate the convergence of the high-fidelity solver 
described in Section \ref{sec:high_fidelity_solver}
for this configuration.
To this end, we consider $512$ realization and take the
average error in the $\IT^0$-norm. the results are
presented in Figure \ref{fig:hferrors16}. 
Based on the errors achieved by the high-fidelity solver, for further results related to this test case we fix $N = 40$, unless otherwise specified. 

\begin{figure}
    \centering    \includegraphics[scale = 0.18]{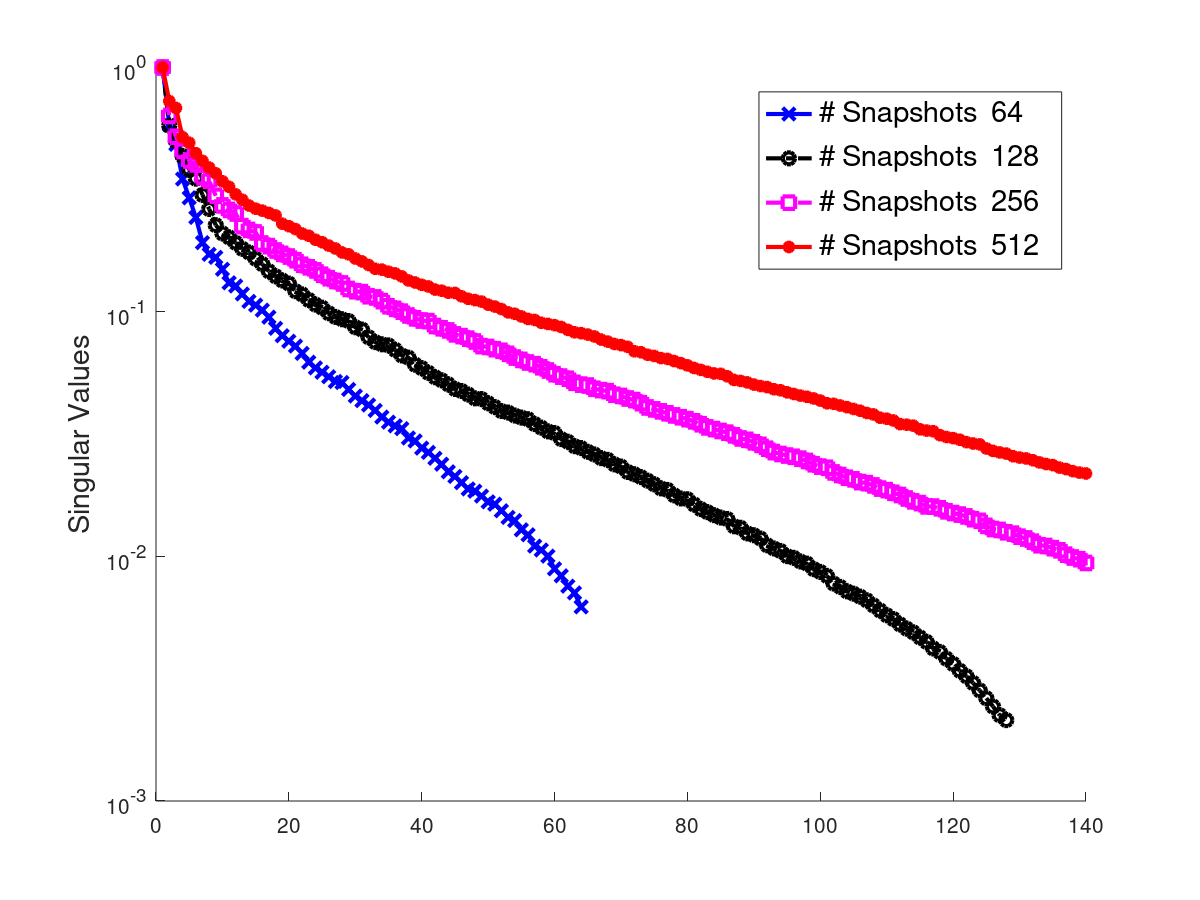}
    \caption{Singular values for the multple arc problem (16 arcs) scaled by the magnitude of the first singular-value and with $N= 70$. }
    \label{fig:singfull16}
\end{figure}

Now we analyze if the multple arc problem is amenable for reduction.
We consider a number of snapshots and inspect the singular-values of the snapshot matrix. 
These are presented in Figure \ref{fig:singfull16}.
We observe that while the singular values decay exponentially,
we are not able to achieve relative small tolerances. Consequently, this configuration is not directly amenable for reduction.

Hence, as pointed out in Section \ref{sec:romconstruction}, we have two key issues.
\begin{enumerate}
    \item The dimension of the perturbation parameter is $12 \times 16 = 192$,
    and many (at least twice the total number of arcs) share the same importance.
    This is the main reason as to why is hard to find a suitable reduced basis for the multiple arc problem. 
    \item 
    Even if we could find good basis for the full problem,
    the cost of constructing this basis can be extremely high,
    as we would require to solve a large number of full-order problems.
  \end{enumerate}

\begin{figure}
    \centering    \includegraphics[scale=0.18]{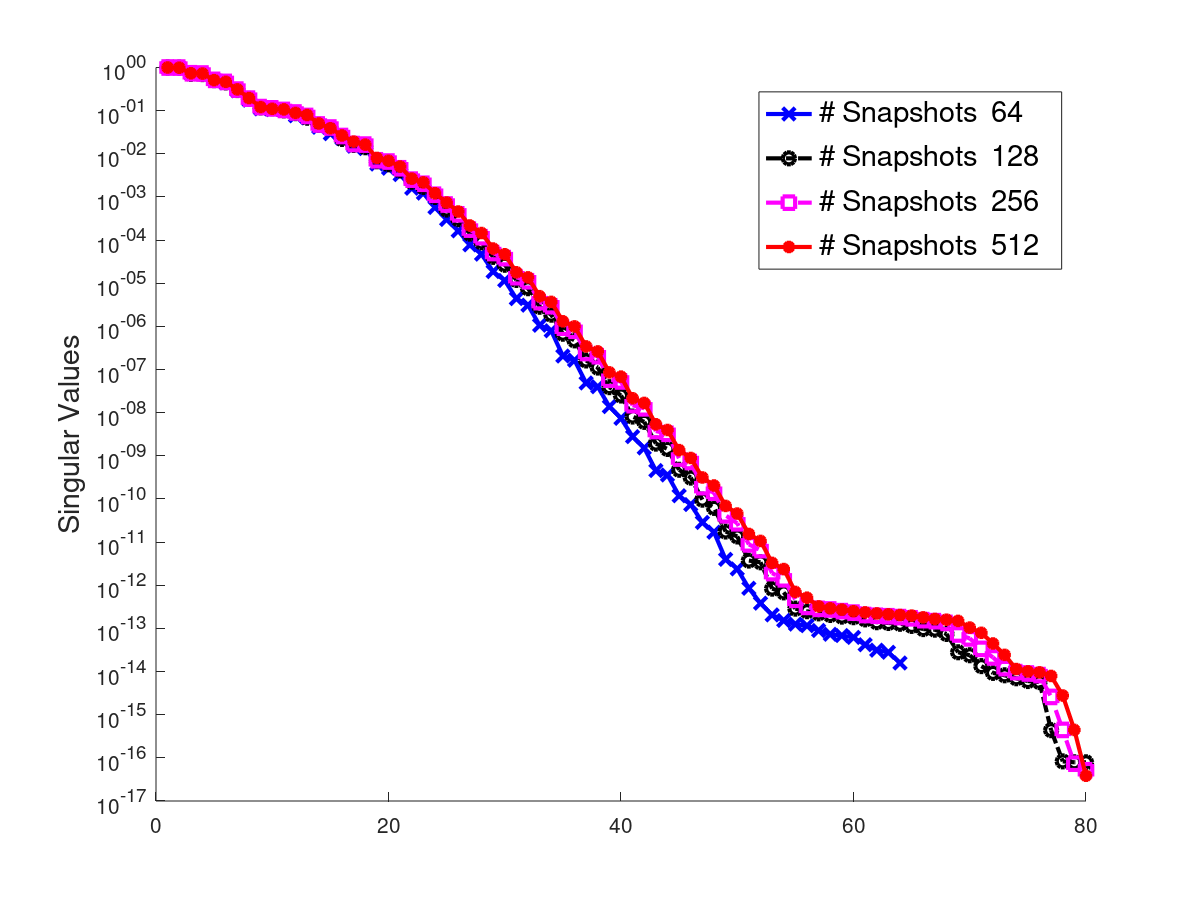}
    \caption{Singular values for a single arc, scaled by the magnitude of the first singular-value, with $N= 100$. }
    \label{fig:singvals1of16}
\end{figure}

Next we study whether the single arc problem is suitable for reduction. 
As explained in Section \ref{sec:eimconstruction}, we consider a parametrization
of the arc that includes the effects of the position, orientation and length as variables
determined by the random parameters, as well as perturbations of the form \eqref{eq:pert16}. 
These results are presented in Figure \ref{fig:singvals1of16}. 
As opposed to the multiple arc problem,
for a single arc we can achieve much smaller relative singular values.  

It remans to dilucidate if the reduced basis for the single arc problem
performs accuratly for the multple arc problem. 
The following results illustrate that this is indeed the case. 

\begin{figure}
  \subfigure[
	Cross-Interactions
  ]{       \includegraphics[width=.45\textwidth]{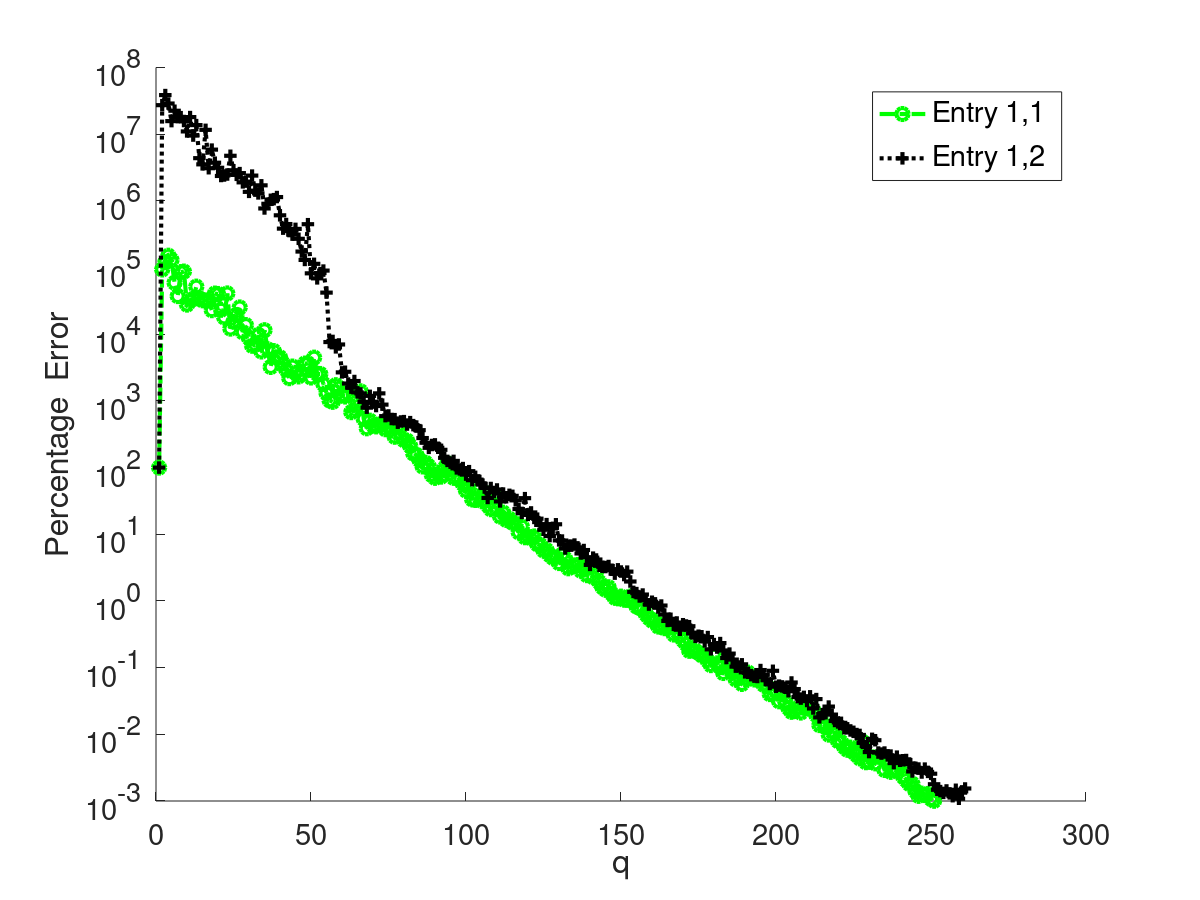}
  }
  \subfigure[
    Self-Interactions
  ]{
  \includegraphics[width=.45\textwidth]{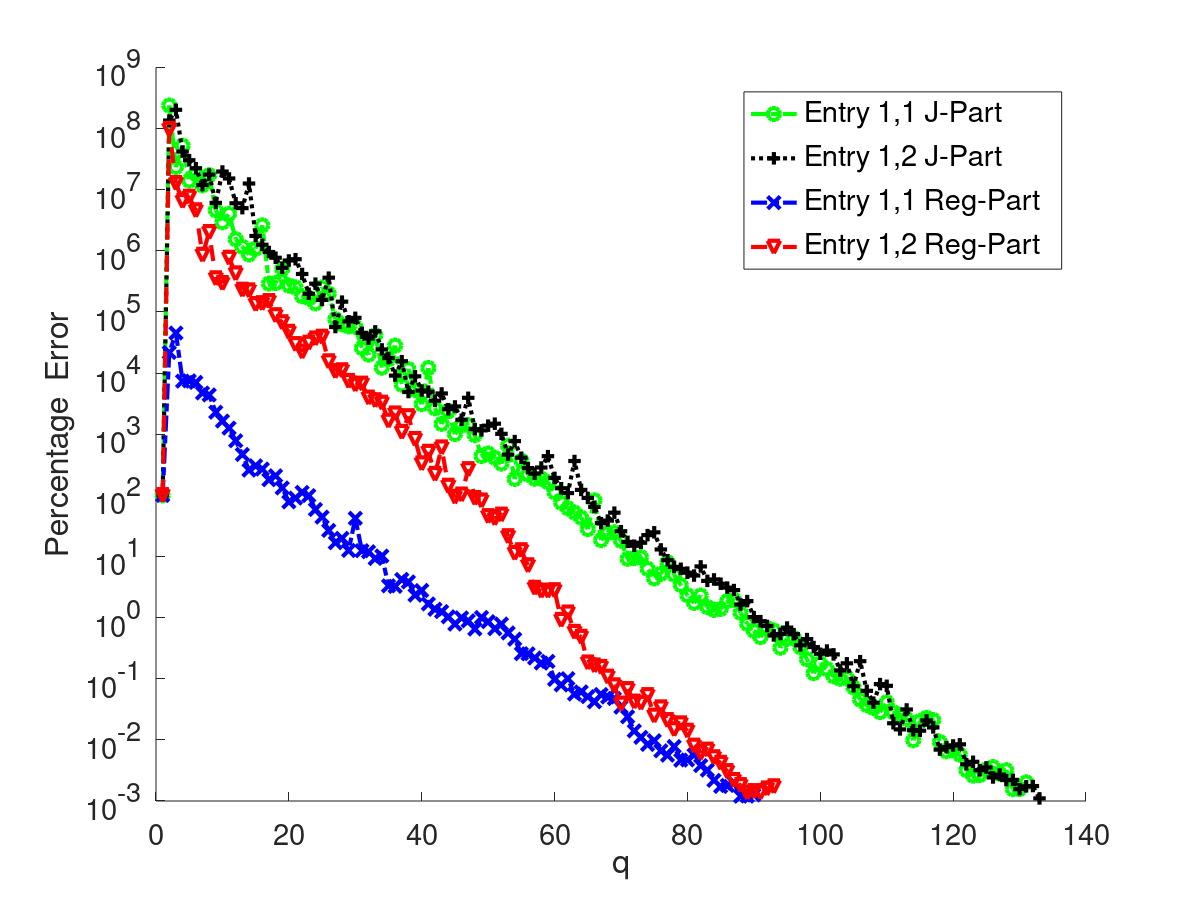}
  }
  \caption{ 
  Percentage error of the interpolation procedure with respect to the number of iterations of Algorithm \ref{euclid}.
}
\label{fig:testeim16}
\end{figure}

Next, we investigate the convergence of the interpolation procedure described in Section \ref{sec:eimconstruction}.
We are interested in the evolution of the relative errors as the number of iterations increases.
Observe that the fundamental solution is a $2 \times 2$ matrix, so we need to compute four scalar interpolations.
However, the pair of diagonal terms exhibits a similar behavior and so does the pair off-diagonal terms.
Thus, we only show results for the entries $(1,1)$ and $(1,2)$.
We also recall that according to Remark \ref{rem:matints}, for the self-interactions
we have to interpolate two type of functions, named the regular-part (for the sake of simplicity referred to as Reg-Part),
and the $J-$part. The corresponding results are presented in Figure \ref{fig:testeim16}.

The remark the following concerning these results: 
\begin{enumerate}
    \item The error appears to decrease exponentially, as the slope of the error curve is approximately constant. 
    \item 
    The cross-interaction matrices need a higher number of terms to achieve the same level of accuracy.
    This should be expected as the cross-interaction maps is of higher dimension as it depends of two parametrically defined arcs
\end{enumerate}

\begin{figure}
    \centering    \includegraphics[scale=0.18]{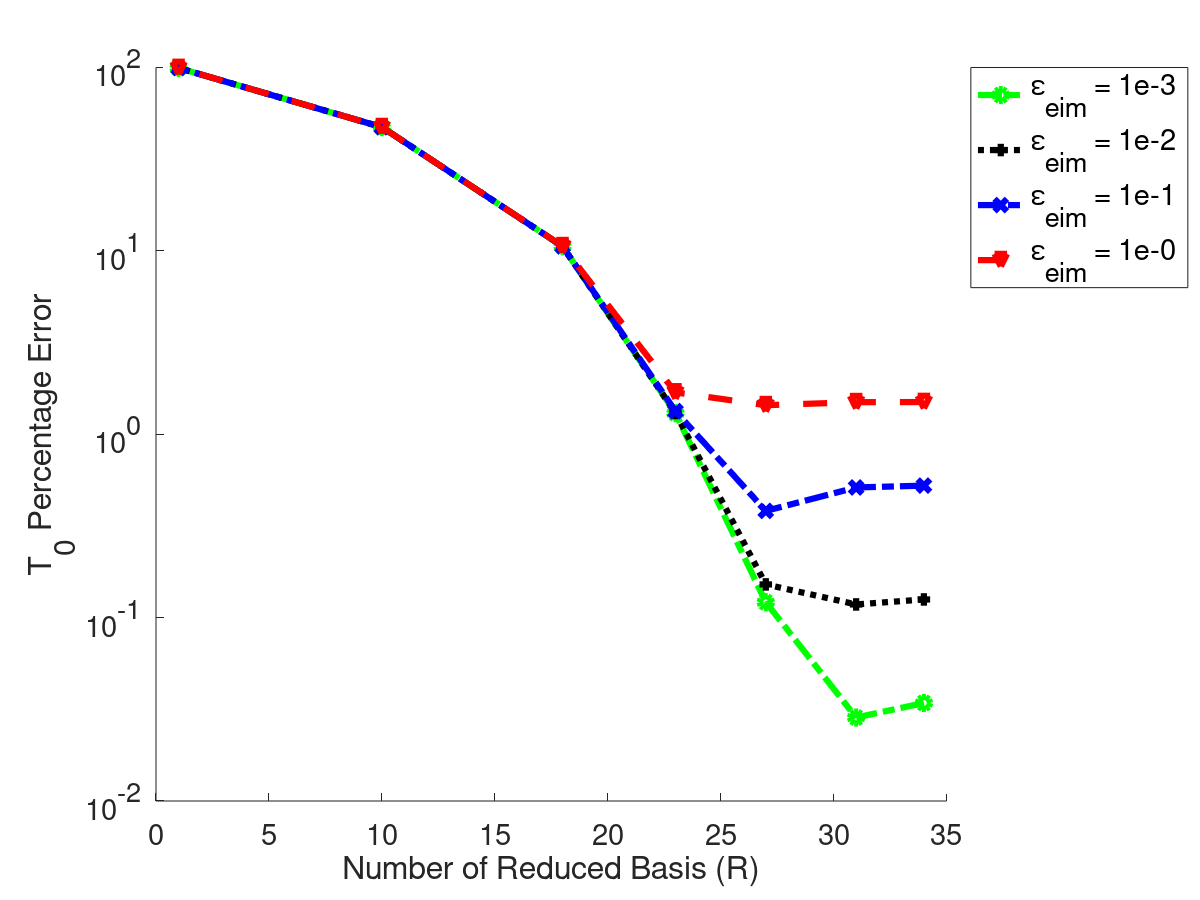}
    \caption{Mean value of the errors of the reduced basis method. }
    \label{fig:test16}
\end{figure}

We present the errors of the reduced bass method for the multiple arc problem.
We recall that according to Sections \ref{sec:pod_rom}
and \ref{sec:eimconstruction} the accuracy of the reduced basis method depends of two parameters,
$\epsilon_\text{svd}$ (cf. Remark, that determines the number of basis, and $\epsilon_\text{eim}$ the tolerance
used to construct the interpolation of the corresponding functions. In Figure \ref{fig:test16} we present
the average error of 56 test cases depending of the number of reduced basis used $R$
(which is determined by $\epsilon_\text{svd}$), and $\epsilon_\text{eim}$.
The solution of the respective linear system is done using the preconditioned GMRES. 

To further illustrate the performance of the method, in Table \ref{tab:test16} we include some extra information regarding these test cases, such as the tolerances utilized,  the computing times: Times-RB (solution time for the reduced basis, Times-HF (solution time for the high-fidelity), and the values of $R$ (number of reduced basis) and $Q$ (number of terms in the interpolation), for the latter we show a weighted mean value, as the value differs for each of the four entries of the fundamental solution, and also if it correspond to a cross or self-interaction. The associated weights are $\frac{M^2-M}{M^2}$, for the cross-interactions, and $\frac{M}{M^2}$, for the self-interactions.

\begin{table}[]
\begin{tabular}{l|l|l|l|l|l|l}
$\epsilon_\text{svd}$ & $R$ & $\epsilon_\text{eim}$ & mean of $Q$ &Time-RB & Time-HF  & Percentage Error \\ \hline
1e-6             & 34  & 1e-3             & 245                & 3.1s              & 6s             & 3.4e-2             \\ \hline
1e-6             & 34  & 1e-1             & 179               & 2s              & 6s             & 5e-1            \\ \hline
1e-3             & 23  & 1e-3             & 245                & 2.8s              & 6s             & 1.3             \\ \hline
1e-3             & 23  & 1e-1             & 179                & 1.8s              & 6s             & 1.3           \\ \hline
1e-1             & 10  & 1e-3             & 245                & 2.5s              & 6s             & 47           \\ \hline
1e-1             & 10  & 1e-1             & 179                & 1.5s              & 6s             & 47          
\end{tabular}
\caption{Convergence results.}
\label{tab:test16}
\end{table}

\subsection{Increasing Number of Arcs}

\begin{figure}
  \subfigure[
	121 Arcs
  ]{       \includegraphics[width=.45\textwidth]{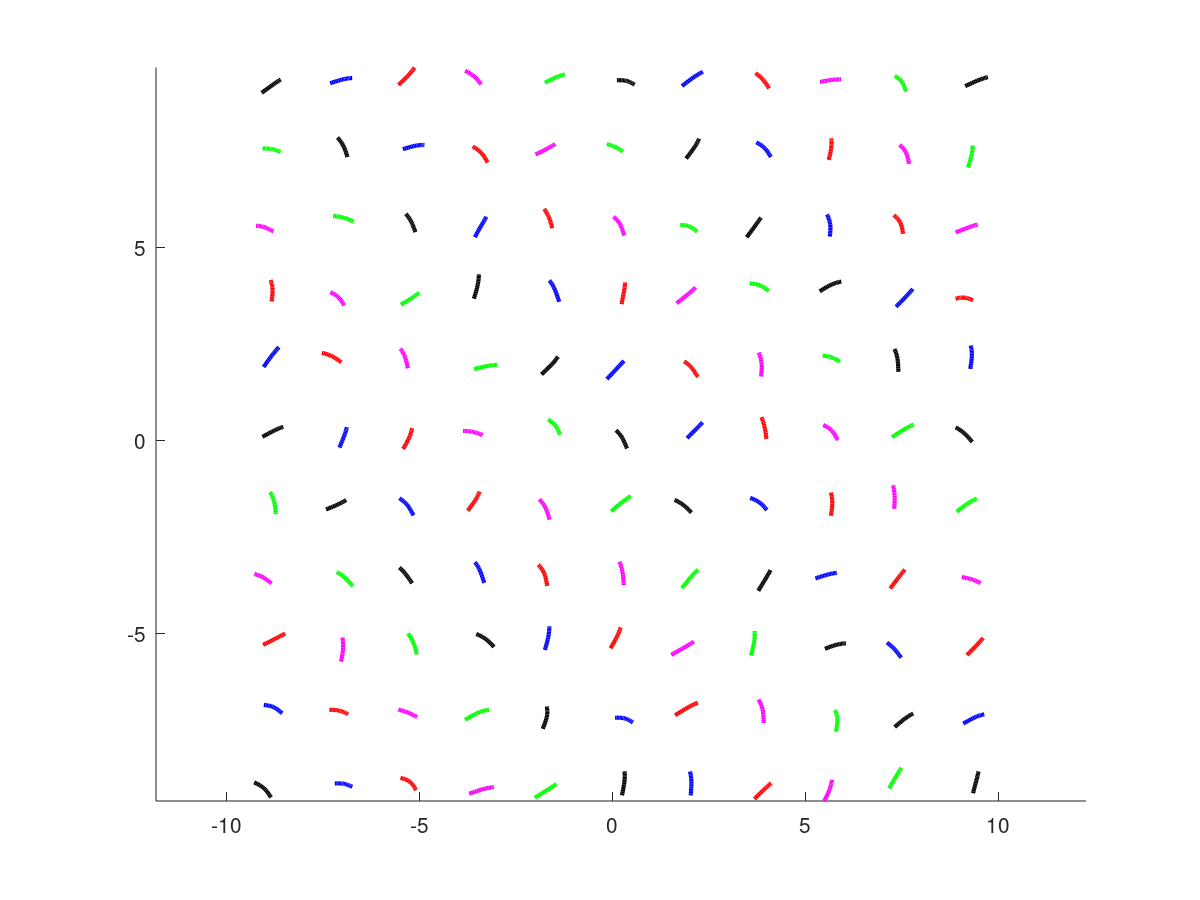}
  }\subfigure[
    256 Arcs
  ]{
  \includegraphics[width=.45\textwidth]{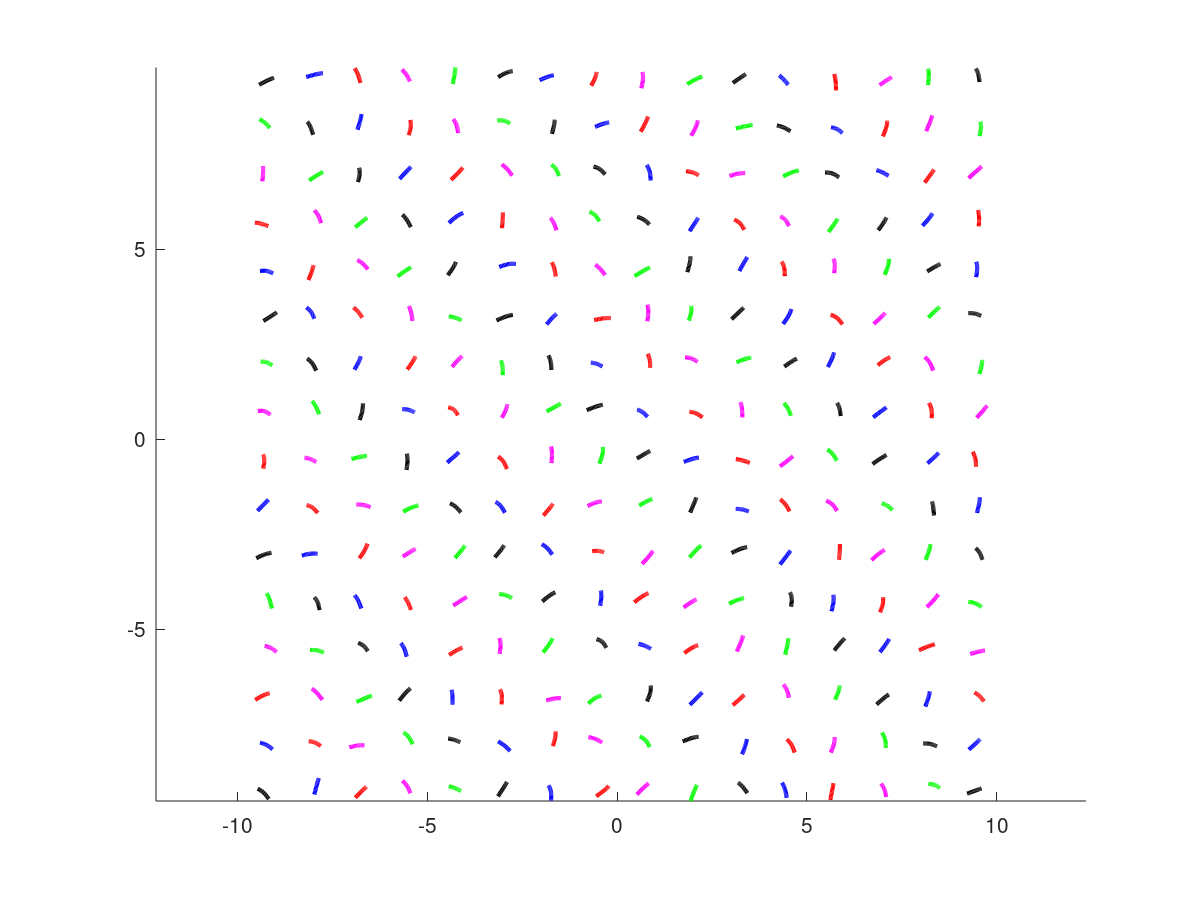}
  }
  \caption{ 
Geometry realizations
}
\label{fig:geomulti}
\end{figure}

We again consider the problem with $\omega = 10,$ $\lambda = 2$, $\mu =1$, and perturbations of the form \eqref{eq:pert16}, but with $c_n = \frac{\min(1, r_\text{max})}{n^{2.5}}$. The arcs will again be positioned in $[-10,10] \times [-10,10],$ but the number of arcs, as well as the global parameters $r_\text{min},r_\text{max}, d_\text{min}, d_\text{max}$ are variables. In particular we will consider different number of arcs, from 36 to 1024, and we adapt the global geometry parameters to ensure that not self crossing occurs, we illustrate some geometry realizations in Figure \ref{fig:geomulti}.

Notice that the size of the arcs decrease as the number of arcs increase, this is somehow equivalent to reduce the frequency as we increase the number of arcs. 

In some of the cases considered in this section (529 and 1024 arcs) it is impossible (given the memory constraints\footnote{All the experiments were performed on a desktop computer I7-4770 with 32gb of RAM.}) to obtain the solutions using the high-fidelity solver, thus we introduce a classical a-posteriori estimation of the error to illustrate the performance of the method. To this end, we define the relative residual as, 
$$
	\text{res} 
	= 
	\sum_{k=1}^M\frac{\| \sum_{j=1}^M \mathbb{A}_{k,j} \mathbb{V}_R^{(\text{rb})} \bm{a}^{\text{(rb)},j}_R - \bm{g}_k \|^2}{\| \bm{g}_k\|^2}
$$
The results are presented in Table \ref{tab:test2}.

\begin{table}[]
\begin{tabular}{l|l|l|l|l|l|l}
M   & R  & mean Q & Time-RB & Time-HF & Percentage Error & Percentage Residual \\ \hline
36  &  22  & 175     & 9s                & 90s                & 0.07 &0.01             \\ \hline
64  &  20  & 140     & 21s                & 290s               & 0.05 & 0.007             \\ \hline
121 &  18  & 112     & 70s               & 479s               & 0.04 & 0.006           \\ \hline
256 & 16 & 88.8     & 311s               & 4640s              & 0.04 & 0.003\\  \hline
529 & 14 & 70 & 1820s & - & - & 0.005 \\  \hline
1024 & 14 & 60 & 11800s & - & - & 0.001
\end{tabular}
\caption{Results for number of arcs (M), with $N= 70$ for the high-fidelity solver, $\epsilon_\text{svd} = \epsilon_\text{eim} = 10^{-3}$.}
\label{tab:test2}
\end{table}

\section{Concluding Remarks}
\label{sec:conclu}
In this work, we present and analyze a reduced basis algorithm for the 
elastic sccattering bby multiple arcs in two space dimensions. The key
insight of the method, which as previously metioned follows \cite{GHS12},
consists first in finding a reduced space for a single shape-parametric, and then use this
as a bulding block for the construction of a reduced space for the multiple open arc.
problem.
Among the advantages of this method, we highlight that once the reduced basis
has been cosntructed and stored in the offline phase, one can use this reduced space
for the computation solution of different multiple arcs configuration, possibly with 
different number of arcs, in the online phase of the reduced basis method. 
Furthermore, based on our previous work \cite{pinto2024shape}, we present a
complete analysis of the method. In particular, in our analysis, we provde an d
argument as to why the reduced basis for the single arc serves for the multiple arc problem.
Even though we have presented a complete description for the multiple arc problem equipped 
with Dirichlet boundary conditions, the exact same analysis can be extended to the same problem equipped
with Neumann boundary conditions. This would lead to a boundary integral formulation characterized
by the presence of hypersingular BIOs. Future work comprises the extension of this work and analysis 
to acoustic and elastic scattering by multiple objects in three dimensions, and the construction of
neural network-based surrogates for the approximation of the corresponding parameter-to-solution 
map as in \cite{hesthaven2018non}.

\bibliography{biblio}
\bibliographystyle{siam}

\end{document}